\documentclass[10pt]{article}
\usepackage[margin=1in]{geometry}

 \usepackage{amsthm}
\usepackage{amsmath}
\usepackage{mathrsfs}
\usepackage{amsfonts}
\usepackage{amssymb}
\usepackage{xcolor}

\newtheorem{lemma}{Lemma}[section]

\newtheorem{theorem}{Theorem}[section]

\newtheorem{definition}{Definition}[section]

\let\Section=\section
\def\section{\setcounter{equation}{0}\Section}
\sf

\begin{document}
\title{Schr\"odinger-Poisson equations with singular potentials in $\mathbb{R}^3$
}
 \author{Yongsheng Jiang  and   Huan-Song Zhou\thanks{Corresponding author.
 \hfil\break\indent
 \hfil\break \indent {\it E-mail address}: hszhou@wipm.ac.cn(H.S.Zhou), jiangys@znufe.edu.cn(Y.S.Jiang)}\\
 \small Wuhan Institute of Physics and Mathematics, Chinese
Academy of Sciences\\
\small  P.O.Box 71010, Wuhan 430071, China}
 \date{}
\maketitle

\noindent {\sc Abstract}:   The existence and $L^{\infty}$ estimate of positive
solutions are discussed for the following Schr\"{o}dinger-Poisson
system
\begin{equation}
\left\{\begin{array}{ll}
 -\Delta u +(\lambda+\frac{1}{|y|^\alpha})u+\phi (x) u =|u|^{p-1}u, \ \ x=(y,z)\in \mathbb{R}^2\times\mathbb{R}, \\
 -\Delta\phi = u^2,\  \lim\limits_{|x|\rightarrow +\infty}\phi(x)=0,  \hfill y=(x_1,x_2) \in
 \mathbb{R}^2 \mbox{ with } |y|=\sqrt{x_1^2+x_2^2},
 \end{array}\right.
\end{equation}
where $\lambda\geqslant0$, $\alpha\in[0,8)$ and
$\max\{2,\frac{2+\alpha}{2}\}<p<5$.



\section{Introduction}

In this paper, we study the following type of Schr\" odinger-Poisson
equations
\begin{equation}\label{eq:1.1}
\left\{\begin{array}{ll}
 -\Delta u +V(x)u+\phi (x) u =|u|^{p-1}u,  \\
 -\Delta\phi = u^2,\ \lim\limits_{|x|\rightarrow +\infty}\phi(x)=0, \ \ \  x=(x_1,x_2,z)\in
 \mathbb{R}^3,
 \end{array}\right.
\end{equation}
where $p\in(2,5)$, and the potential function $V(x)$ is of the form
\begin{description}
 \item $(V)$ $V_{\lambda}(x)=\lambda+\frac{1}{|y|^\alpha}$,\ $\lambda\geqslant0$,\ $\alpha\in[0,8)$ \
  and $|y|=\sqrt{x_1^2+x_2^2}$.
\end{description}
Problem (\ref{eq:1.1}) can be viewed as the stationary problem of
 the following coupled Schr\"odinger-Poisson system:
\begin{equation}\label{eq:1.1.1}
\left\{\begin{array}{ll}
 i\psi_{t}-\Delta \psi +\phi (x) \psi =f(|\psi|)\psi,  \\
 -\Delta\phi = |\psi|^2,\ \lim\limits_{|x|\rightarrow +\infty}\phi(x)=0, \ \ \  x\in
 \mathbb{R}^3,
 \end{array}\right.
\end{equation}
where $f(|\psi|)\psi=|\psi|^{p-1}\psi+\omega_0\psi$, $\omega_0>0$,
$2<p<5$ and
$\psi:\mathbb{R}^3\times\mathbb{R}\longrightarrow\mathbb{C}$. In fact, motivated by
\cite{VBFD}, we may seek a solution of (\ref{eq:1.1.1}) with the following type:
$$\psi(x,t)=u(x)e^{i(\eta(x)+\omega t)},\ \ u(x)\geq0,\ \eta(x)\in \mathbb{R}/2\pi\mathbb{Z},\ \omega\geq\omega_0.$$
Then, by (\ref{eq:1.1.1}),  $u$ should satisfy a system
\begin{equation}\nonumber
\left\{\begin{array}{ll}
 -\Delta u +(\omega-\omega_0+|\nabla\eta(x)|^2)u+\phi (x) u =|u|^{p-1}u,  \\
 u\Delta\eta(x)+2\nabla u\nabla\eta=0,\\
 -\Delta\phi = u^2,\ \lim\limits_{|x|\rightarrow +\infty}\phi(x)=0, \ \ \  x\in
 \mathbb{R}^3.
 \end{array}\right.
\end{equation}
Furthermore, similar to \cite{MBadVBenSRol-reprent,BFD1}, for $x\in \mathbb{R}^3$, if we let
$u(x)=u(y,z)=u(|y|,z)$ and
\begin{equation}\nonumber
\eta(x)=\left\{\begin{array}{ll}
 \arctan(x_2/x_1),\ \  \text{ if }x_1>0, \\
\arctan(x_2/x_1)+\pi,\ \  \text{ if }x_1<0, \\
\pi/2,\ \  \text{ if }x_1=0\  \text{ and } x_2>0,\\
-\pi/2,\ \  \text{ if }x_1=0\  \text{ and } x_2<0,\\
 \end{array}\right.
\end{equation}
it is easy to see that $\eta(x)\in C^2(\mathbb{R}^3\setminus T_-)$,
where $T_-:=\{(x_1,x_2,z)\in\mathbb{R}^3:x_1=0,x_2\leq0\}$. By a
simple calculation we know that
$$\Delta\eta(x)=0,\ \ \nabla\eta(x)\cdot\nabla u(x)=0,\ \ |\nabla\eta(x)|=\frac{1}{|y|^2},\ \text{ for    } x\in\mathbb{R}^3\setminus T_- .$$
These show that $u(|y|,z)$ is actually a nonnegative solution of
(\ref{eq:1.1}) with $\alpha=2$ and $\lambda=\omega-\omega_0$. Furthermore,
$\psi(x)$ solves (\ref{eq:1.1.1}) with angular momentum:
\begin{equation}\nonumber
M(\psi)=\text{Re}\int_{\mathbb{R}^3}i\bar{\psi}x\wedge\nabla\psi
dx=-\int_{\mathbb{R}^3}u^2x\wedge\nabla v(x) dx=-(0,0,|u|_{L^2}^2).
\end{equation}

For problem (\ref{eq:1.1}), more and more results have been published under various conditions on the potential
function $V(x)$ and on the nonlinear term $|u|^{p-1}u$, for examples,
if $V(x)=const$, that is $\alpha=0$ in $(V)$, the non-existence of
nontrivial solution of (\ref{eq:1.1}) for $p\not\in(1,5)$ was proved
in \cite{TDApriDMugn-AdvNS} by a Pohozaev type identity, a radially
symmetric positive solution  was obtained
in \cite{GMColi-CAA} and \cite{TDApriDMugn-ProcRSES} for $p\in[3,5)$, etc.
It is known that we may find a nontrivial weak solution of problem (\ref{eq:1.1}) by looking for a nonzero critical point of the related variational functional of problem (\ref{eq:1.1}).  It is also
known that the weak limit of a so-called Palais-Smale sequence ((PS) sequence, in short) of the variational functional  is usually a
weak solution, but it may be a trivial
solution unless we can prove that the variational functional satisfies the Palais-Smale condition ((PS)
condition, in short), that is, a (PS) sequence has a strongly convergent
subsequence. However, without condition
(\ref{eq:1.2}) below, it seems very difficult to show a (PS)
sequence converges strongly. In this paper, instead of trying to
prove the (PS) condition,  we adapt a trick used in
\cite{BadMGuiMRolS-AdvDE}, which is essentially a version of the
concentration-compactness principle due to \cite{SSol-AnnIHPA}, to
show directly that the weak limit of a (PS) sequence is indeed a
nontrivial solution. For this purpose, we have to ensure that the
(PS) sequence obtained by the deformation Lemma \cite{MWil} is
nonnegative and {$\phi(x)$ is bounded in $D^{1,2}(\mathbb{R}^3)$}, this is because there
is a term $\phi(x)u$ appearing in problem (\ref{eq:1.1}), which is usually
called  a nonlocal term. As a by-product, in this paper we provide a
simple approach for getting a nonnegative (PS) sequence and a bound of $\phi(x)$
 in $D^{1,2}(\mathbb{R}^3)$, see Lemma \ref{L:2.6}, this may be useful
in certain situations. Note that in
\cite{BadMBenciVRolS-JEMS,BV,DT,BadMGuiMRolS-AdvDE} the authors
studied the single stationary Schr\"{o}dinger equation, that is, the
first equation of (\ref{eq:1.1}) with $\phi(x)=0$ (i.e. without
nonlocal term), in this case it is not necessary to seek a
nonnegative (PS) sequence, see e.g.
\cite{BadMBenciVRolS-JEMS,BadMGuiMRolS-AdvDE}. It seems no any
results for Schr\"oding-Poisson system (\ref{eq:1.1}) under
condition ($V$) with $\alpha>0$. We should mention that our results
of this paper cover the case of $\alpha=0$, that is, the constant
potential case. In this paper, we give also a priori estimate for
 solutions of (\ref{eq:1.1}), see Lemma \ref{Le:4.4.0}, and get also a classical solution (except $|y|=0$) for
(\ref{eq:1.1}) with
$\lambda=0$, $\alpha\in(0,8)$ and $\max\{2,\frac{2+\alpha}{2}\}<p<5$.

For problem (\ref{eq:1.1}) with constant potential, i.e. taking
$\alpha=0$ in $(V)$, the existence and nonexistence results were
established by Ruiz in \cite{DRuiz-JFA}, he proved that
(\ref{eq:1.1}) has always a positive radial solution if $p\in(2,5)$
and does not admit any nontrivial solution if $p\leqslant2$. A
ground state for (\ref{eq:1.1}) with $p\in(2,5)$ was proved in
\cite{AAzzoAPom-Jmaa}. The existence of non-radially symmetric
solution was shown in \cite{PDAven-AdvNS} and  multiple solutions
for (\ref{eq:1.1}) were obtained in
\cite{AAmbroDRuiz-CCM,GMColi-CAA}. If the potential $V(x)$ is not a
constant, problem (\ref{eq:1.1}) has been studied in
\cite{AAzzoAPom-Jmaa} for $p\in (3,5)$ and \cite{LgZhaoFkZhao-Jmaa}
for $p \in (2,3]$.  (\ref{eq:1.1}) with more general
nonlinearities has been studied in
\cite{ASS,AAP,Mug,ZPWHSZhou-DisCDS,YZD}, etc. To ensure that the
variational functional of problem (\ref{eq:1.1}) satisfies the (PS)
condition, in the papers
\cite{AAzzoAPom-Jmaa},\cite{LgZhaoFkZhao-Jmaa} the following
conditions are assumed
\begin{equation} \label{eq:1.2}
V(x)\leqslant V_{\infty}=\liminf\limits_{|x|\rightarrow\infty}V(x),
 \end{equation}
 \begin{equation} \label{eq:1.2.0}
2V(x)+(\nabla V(x),x) \geqslant0 \text{ a.e. }
 x\in\mathbb{R}^3.
 \end{equation}
It is clear that the above conditions are not true for the potential
given by $(V)$. So, we cannot follow the same tricks as that of
\cite{AAzzoAPom-Jmaa,LgZhaoFkZhao-Jmaa} to deal with problem
(\ref{eq:1.1}). Without condition (\ref{eq:1.2.0}), it seems
difficult even in showing that a (PS) sequence is bounded in the working
Sobolev space, specially in the case of $p\in(2,3)$. Motivated by
\cite{MBadVBenSRol-reprent}, here we try to find a bounded and
nonnegative (PS) sequence directly from the well-known
deformation Lemma (\cite{MWil}, Lemma 2.3).\\

Before stating our main results, we introduce some notations,
definitions and recall some properties of the solution of the second
equation (Poisson equation) in (\ref{eq:1.1}). For
$\alpha\geq0$ and $x=(y,z)\in\mathbb{R}^2\times\mathbb{R}$, define

\begin{equation}\label{eq:1.4}
E=\{u(x)\in \textrm{D}^{1,2}(\mathbb{R}^3):u(x)=u(|y|,z)\text{ and
}\int_{\mathbb{R}^3}\frac{u^{2}}{|y|^{\alpha}}dx<\infty
\},
\end{equation}
and for $\lambda>0$, we denote

\begin{equation}\nonumber
H=\{u\in E:\lambda\int_{\mathbb{R}^3}u^{2}dx<\infty \}.
\end{equation}
Clearly $H\subset E$, $H\subset H^1(\mathbb{R}^3)$ and $H$ is a Hilbert
space, its scalar product and norm are given by
\begin{equation}\label{eq:1.6}
\langle u, v\rangle_{H}=\int_{\mathbb{R}^3}[\nabla u\nabla
v+V_{\lambda}(x)uv]dx\ \text{ and } \ \|u\|_{H}^2=\langle u,
u\rangle_{H},
\end{equation}
respectively, where $V_{\lambda}(x)=\lambda+\frac{1}{|y|^{\alpha}}$. \\

Throughout this paper, we denote the standard norms of
$H^1(\mathbb{R}^3)$ and $L^p(\mathbb{R}^3)$ ($1\leqslant
p\leqslant+\infty$) by $\|\cdot\|$ and $\|\cdot\|_p$, respectively.
Then, (\ref{eq:1.6}) implies that $\|\cdot\|_{H}$ is an equivalent
norm of $\|\cdot\|$ if $\alpha=0$.\\

 By Lemma 2.1 of \cite{DRuiz-JFA}, we know that
$-\Delta\phi(x)=u^2$ has a unique solution in
$D^{1,2}(\mathbb{R}^3)$ with the form of
\begin{equation}\label{eq:1.7}
\phi(x):=\phi_{u}(x)=\int_{\mathbb{R}^3}\frac{u^2(y)}{|x-y|}dy,\text{
for any } u\in L^{\frac{12}{5}}(\mathbb{R}^3),
\end{equation}
and
\begin{equation}\label{eq:1.8}
|\nabla\phi_{u}(x)|_2\leqslant C|u|_{12/5}^2,\ \
\int_{\mathbb{R}^3}\phi_{u}(x)u^2dy\leqslant C|u|^4_{12/5}.
\end{equation}
For $\lambda>0$ and $u\in H$, we can define
 the variational functional of problem (\ref{eq:1.1}) as follows:
\begin{equation} \label{eq:1.10}
I(u):=I_\lambda(u)=\frac{1}{2}\int_{\mathbb{R}^{3}}({|\nabla
u|}^{2}+V_{\lambda}(x)u^2)dx+\frac{1}{4}\int_{\mathbb{R}^{3}}\phi_u(x)u^2dx
-\frac{1}{p+1}\int_{\mathbb{R}^{3}}{|u|}^{p+1}dx.
 \end{equation}
Since (\ref{eq:1.8}), $I_\lambda$ is well defined on $H$ and
$I_\lambda\in C^1(H,\mathbb{R})$ with
\begin{equation}\label{eq:1.11}
(I'_\lambda(u),v)=\int_{\mathbb{R}^{3}}(\nabla u\nabla
v+V_{\lambda}(x)uv)dx+\int_{\mathbb{R}^{3}}\phi_u(x)uvdx
-\int_{\mathbb{R}^{3}}{|u|}^{p-1}uvdx
 \end{equation}
for all $v\in H$ with $\lambda>0$ and $p\in (1,5)$.
Furthermore, {it is known that a weak
solution of (\ref{eq:1.1}) corresponds to a nonzero critical point
of the functional $I$ in $H$} if $\lambda>0$.

However, if $\lambda=0$, then $H=E$. In this case,
 {(\ref{eq:1.7}) (\ref{eq:1.8}) are not always true for $u\in E$. Therefore, the integrations
 $\int_{\mathbb{R}^{3}}{|u|}^pdx$,
$\int_{\mathbb{R}^{3}}\phi_u(x)u^2dx$ and
$\int_{\mathbb{R}^{3}}\phi_u(x)uvdx$ may not be well defined for
$u$, $v\in E$}.

In this paper, we want to establish some existence results for problem (\ref{eq:1.1}) for both $\lambda >0$ and $\lambda =0$.
To this end, we set
\begin{equation}\label{T}
T=\{x\in \mathbb{R}^3:|y|=0\} \text{ where }|y|=\sqrt{x_1^2+x_2^2}.
\end{equation}
Hence, by an approximation procedure, see Section 4, we can find a weak solution $u\in E$ of
(\ref{eq:1.1}) with $\lambda=0$  in the sense of
\begin{equation} \label{eq:1.9.1}
\int_{\mathbb{R}^{3}}(\nabla u\nabla
\varphi+\frac{1}{|y|^\alpha}u\varphi)dx+\int_{\mathbb{R}^{3}}\phi_u(x)u\varphi
dx =\int_{\mathbb{R}^{3}}{|u|}^{p-1}u\varphi dx,\ \text{ for }
\varphi\in C^\infty_0(\mathbb{R}^3\setminus T).
 \end{equation}
Note that $\int_{\mathbb{R}^3} \frac{1}{|y|^\alpha} u\varphi dx$ may be not  integrable for $u\in E$ and $\varphi \in
C^\infty_0(\mathbb{R}^3)$, this is why we take $\varphi \in C^\infty_0(\mathbb{R}^3\setminus T)$ above instead of $\varphi \in C^\infty_0(\mathbb{R}^3)$.
So, it is reasonable for us to define a weak solution for  (\ref{eq:1.1}) as follows.

\begin{definition}\label{def:1.1}
$u\in E\setminus\{0\}$ is said to be a weak solution of
(\ref{eq:1.1}) with $\lambda\geqslant0$ if $\phi_{u}\in
D^{1,2}(\mathbb{R}^3)$ and $u$ satisfies
\begin{equation} \label{eq:1.9}
\int_{\mathbb{R}^{3}}[\nabla u\nabla
\varphi+(\frac{1}{|y|^{\alpha}}+\lambda)u\varphi]dx+\int_{\mathbb{R}^{3}}\phi_u(x)u\varphi
dx =\int_{\mathbb{R}^{3}}{|u|}^{p-1}u\varphi dx
 \end{equation}
 for all $\varphi\in C^\infty_0(\mathbb{R}^3\setminus T)$.
\end{definition}

We mention that the above definition also enables us to get a
classical solution. In fact, if $u\in E$ and $ \phi_u\in
D^{1,2}(\mathbb{R}^3)$ satisfies (\ref{eq:1.9}),
 {by using our Lemmas \ref{L:4.1} and
\ref{Le:4.3}}, as well as Theorems 8.10 and 9.19 in
\cite{DavidGilbarg-Neil.S.Trudinger}, we can prove that $u\in
C^2(\mathbb{R}^3\setminus T)$, that is, $u$ is a classical solution of
(\ref{eq:1.1}), see Theorem \ref{th:3.1} in section 3.

For the following single Schr\"{o}dinger equation
\begin{equation} \label{eq:1.12}
-\Delta u +\frac{u}{|y|^\alpha} =f(u),
\,\,\,x=(x_1,x_2,\cdots,x_N)\in \mathbb{R}^N,\,N\geqslant3 \
 \end{equation}
with $|y|=\sqrt{\Sigma_{k=1}^{N+1-i}x^2_k}$, $i<N$, the authors of
paper \cite{BadMBenciVRolS-JEMS} proved that (\ref{eq:1.12}) has a
nontrivial solution in $H^1(\mathbb{R}^N)$ if $\alpha=2$,
$N>i\geqslant2$ and $f(t)$ is supposed to have some kinds of double
powers behavior which ensure that $F(u)=\int_0^{u}f(s)ds$ is
well defined in $L^1(\mathbb{R}^N)$ when $u\in
D^{1,2}(\mathbb{R}^N)$. In \cite{BadMBenciVRolS-JEMS}, the authors used a variational method to seek
first a nontrivial solution of
(\ref{eq:1.12}) in $D^{1,2}(\mathbb{R}^N)$, then they proved this solution is in
$L^2(\mathbb{R}^N)$. Formally, (\ref{eq:1.12}) is nothing but the first
equation of problem (\ref{eq:1.1}) by taking $\lambda=0$, $N=3$ and
getting rid of the nonlocal term $\phi(x)u$. However, even for $f(u)=|u|^{p-1}u$
with $p\in(2,5)$, $F(u)$ is not well defined in $D^{1,2}(\mathbb{R}^N)$,
then the method and results of \cite{BadMBenciVRolS-JEMS} do not
work for our problem. For these reasons, it seems difficult to choose a
working space to solve (\ref{eq:1.1}) directly if $\lambda=0$. In this
paper, we prove first that (\ref{eq:1.1}) has always a solution
$u_\lambda$ in $H^1(\mathbb{R}^3)$ for each $\lambda>0$, then show
that $\{u_{\lambda}\}$ (as a sequence of $\lambda$) is bounded in
$E$,
as mention above we can finally use an approximation process
to get a weak solution of (\ref{eq:1.1}) for $\lambda=0$ in the
sense of (\ref{eq:1.9.1}).
\\

The main results of this paper can be stated now as follows:
\begin{theorem}\label{th1.1}
Let $\alpha\in[0,8)$, $\max\{2,\frac{2+\alpha}{2}\}<p<5$ and
condition $(V)$ be satisfied. Then, problem (\ref{eq:1.1}) has at least a
positive solution $u_\lambda\in H \cap
C^2(\mathbb{R}^3\setminus T)$ for every $\lambda>0$.
Furthermore, if $\lambda\in(0,1]$, there exists $C>0$ which is independent of
$\lambda\in(0,1]$ such that the solution $u_\lambda$ satisfies
$$\|\nabla u_\lambda\|^2_{2}+\int_{\mathbb{R}^{3}}\phi_{u_\lambda}u_\lambda^2 dx<C.$$
\end{theorem}

\begin{theorem}\label{th1.2} For $\lambda=0$, let $\alpha\in[0,8)$ and
$\max\{2,\frac{2+\alpha}{2}\}<p<5$. Then, problem (\ref{eq:1.1}) has at
least a positive solution $u\in E\cap
C^2(\mathbb{R}^3\setminus T)$ in the sense of (\ref{eq:1.9.1}).
\end{theorem}

\section{Bounded nonnegative (PS) sequence}
In this section, $\lambda>0$ is always assumed. Our aim is to known
how the functional $I_\lambda$ defined in (\ref{eq:1.10}) has always
a bounded nonnegative (PS) sequence at some level $c>0$ in $H$. As
 mentioned in the introduction, the authors in
\cite{MBadVBenSRol-reprent} developed an approach to get a bounded
(PS) sequence for the single equation (\ref{eq:1.12}) with certain
nonlinearities. By improving some techniques used in
\cite{MBadVBenSRol-reprent}, we are able to obtain a bounded
nonnegative (PS) sequence for (\ref{eq:1.1}), the nonnegativity of the
(PS) sequence helps us to estimate the related term caused by the
nonlocal term $\phi(x)u$, which leads to a nonzero weak limit of
the (PS) sequence. Let us recall first a deformation lemma from
\cite{MWil}.

\begin{lemma}\label{L:2.1} (\cite{MWil},Lemma 2.3)
Let $X$ be a Banach space, $\varphi\in C^1(X,\mathbb{R})$, $S\subset
X$, $c\in\mathbb{R}$, $\varepsilon,\delta>0$ such that for any
$u\in\varphi^{-1}([c-2\varepsilon,c+2\varepsilon])\cap
S_{2\delta}$\text{:} $\varphi{'}(u)\geqslant8\varepsilon/\delta$.
Then there exists $\eta\in C([0,1]\times X,X)$ such that
\begin{description}
 \item $\rm (i)$  $\eta(t,u)=u$, if $t=0$ or $u\notin\varphi^{-1}([c-2\varepsilon,c+2\varepsilon])\cap
S_{2\delta}$.
\item $\rm (ii)$  $\eta(1,\varphi^{c+\varepsilon}\cap
S)\subset\varphi^{c-\varepsilon}$, where
$\varphi^{c\pm\varepsilon}=\{u\in X:\varphi(u)\leq
c\pm\varepsilon\}$.
\item $\rm (iii)$ $\eta(t,\cdot)$ is an homeomorphism of $X$, for any
$t\in[0,1]$.
\item $\rm (iv)$  $\varphi(\eta(\cdot,u))$ is non increasing, for any
$u\in X$.
\end{description}
\end{lemma}
Now, we give some lemmas, by which Lemma \ref{L:2.1} can be used to
get a desirable (PS) sequence.
\begin{lemma}\label{L:2.2} Let $M>0$ be a constant. If
$u_1$,$u_2\in H$ with $\lambda>0$ and $\|u_1\|_{H}$,
$\|u_2\|_{H}\leqslant M$, then there exist $C:=C(M,p)>0$ such that
\begin{equation} \label{eq:2.1}
\|I{'}(u_1)-I{'}(u_2)\|_{H{'}}\leqslant
C\left(\|u_1-u_2\|_{H}+\|u_1-u_2\|^{3}_{H}\right).
 \end{equation}
\end{lemma}
{\it Proof.} By (\ref{eq:1.11}) and (\ref{eq:1.6}),
\begin{equation*}
\begin{split}
 \left\langle I{'}(u_1)-I{'}(u_2),\psi\right\rangle_{H}&=\langle
u_1-u_2,\psi\rangle_{H}+\int_{\mathbb{R}^{3}}(\phi_{u_1}u_1-\phi_{u_2}u_2)\psi
dx\nonumber\\
&-\int_{\mathbb{R}^{3}}({|u_1|}^{p-1}u_1-{|u_1|}^{p-1}u_1)\psi dx,
\end{split}
\end{equation*}

hence (\ref{eq:2.1}) is proved if we have that
\begin{eqnarray}\label{eq:2.2}
\left|\int_{\mathbb{R}^{3}}({|u_1|}^{p-1}u_1-{|u_2|}^{p-1}u_2)\psi
dx\right|\leqslant C\|u_1-u_2\|_{H}\|\psi\|_{H},
\end{eqnarray}
\begin{eqnarray}\label{eq:2.3}
\int_{\mathbb{R}^{3}}(\phi_{u_1}u_1-\phi_{u_2}u_2)\psi dx\leqslant
C\left(\|u_1-u_2\|_{H}+\|u_1-u_2\|^{3}_{H}\right)\|\psi\|_{H}.
\end{eqnarray}
Indeed, using Taylor's formula and H\"{o}lder inequality as well as
Minkovski inequality, we see that there is a function $\theta$ with
$0<\theta<1$ such that
\begin{eqnarray}
&&\left|\int_{\mathbb{R}^{3}}({|u_1|}^{p-1}u_1-{|u_2|}^{p-1}u_2)\psi
dx\right|\leqslant p\|u_1-u_2\|_{p+1}\|\psi\|_{p+1}\|(\theta
u_1+(1-\theta)u_2)\|^{p+1}_{p+1}\nonumber\\
&\leqslant&p(
\|u_1\|_{p+1}+\|u_2\|_{p+1})^{p+1}\|u_1-u_2\|_{p+1}\|\psi\|_{p+1}
\leqslant p (2M)^{p+1}\|u_1-u_2\|_{p+1}\|\psi\|_{p+1}\nonumber,
\end{eqnarray}
hence (\ref{eq:2.2}) is obtained. To prove (\ref{eq:2.3}), we let
$\upsilon=u_2-u_1$, it follows from (\ref{eq:1.7}) that
\begin{align}\label{eq2-4}
\int_{\mathbb{R}^{3}}(\phi_{u_2}u_2-\phi_{u_1}u_1)\psi
dx=J_1+J_2+J_3+J_4+J_5,
\end{align}
where
\begin{eqnarray*}
J_1=\int_{\mathbb{R}^{3}}\frac{\upsilon^2(y)\upsilon(x)\psi(x)}{|x-y|}dx
dy\leqslant C\|\upsilon\|_{H}^3\|\psi\|_{H}
\end{eqnarray*}
\begin{eqnarray*}
J_2=\int_{\mathbb{R}^{3}}\frac{\upsilon^2(y)u_1(x)\psi(x)}{|x-y|}dx
dy\leqslant C\|\upsilon\|_{H}^2\|u_1\|_{H}\|\psi\|_{H}
\end{eqnarray*}
\begin{eqnarray*}
J_3=\int_{\mathbb{R}^{3}}\frac{u_1^2(y)\upsilon(x)\psi(x)}{|x-y|}dx
dy\leqslant C\|u_1\|_{H}^2\|\upsilon\|_{H}\|\psi\|_{H}
\end{eqnarray*}

\begin{eqnarray*}
J_4=2\int_{\mathbb{R}^{3}}\frac{u_1(y)u_1(x)\upsilon(y)\psi(x)}{|x-y|}dx
dy\leqslant C\|u_1\|_{H}^2\|\upsilon\|_{H}\|\psi\|_{H}
\end{eqnarray*}

\begin{eqnarray*}
J_5=2\int_{\mathbb{R}^{3}}\frac{u_1(y)\upsilon(y)\upsilon(x)\psi(x)}{|x-y|}dx
dy\leqslant C\|u_1\|_{H}\|\upsilon\|_{H}^2\|\psi\|_{H}
\end{eqnarray*}
here we used the following Hardy-Littlewood-Sobolev inequality
(Theorem 4.3 of \cite{EHLiebMLoss-analysis})
\begin{eqnarray*}
\left|\int_{\mathbb{R}^{N}}\int_{\mathbb{R}^{N}}f(x)|x-y|^{-\lambda}h(y)dxdy\right|\leqslant
C(n,\lambda,p)\|f\|_p\|h\|_r,
\end{eqnarray*}
where $p,r>1$ and $0<\lambda<N$ with
$\frac{1}{p}+\frac{\lambda}{N}+\frac{1}{r}=2$, $f\in
L^p(\mathbb{R}^N)$ and $h\in L^r(\mathbb{R}^N)$, the sharp constant
$C(N,\lambda,p)$, independent of $f$ and $h$. Then (\ref{eq:2.3})
holds by (\ref{eq2-4}) and the estimates for $J_1$ to $J_5$. Hence
Lemma \ref{L:2.2} is proved. $\Box$\\

Before giving our next lemma, we recall some basic properties of
$\phi_u(x)$ given by (\ref{eq:1.7}). Let
$$\ u_t:=u_t(x)=t^2u(tx)\ \text{for}\ t>0\ \text{and}\
 x\in \mathbb{R}^3,$$
then $\ u(x)=(u_t)_\frac{1}{t}(x)=(u_\frac{1}{t})_t(x)$ and
\begin{equation}\label{eq:2.3.1}
\|\nabla u_t\|^2_{2}=t^3\|\nabla u\|^2_{2},\quad\quad
\|u_t\|^p_{p}=t^{2p-3}\|u\|^p_{p}\ \text{for}\ 1\leqslant p<\infty,
\end{equation}

\begin{equation}\label{eq:2.3.2}
\int_{\mathbb{R}^{3}}\phi_{u_t}u_{t}^{2}dx=t^3\int_{\mathbb{R}^{3}}\phi_{u}u^{2}dx,
\quad\quad\int_{\mathbb{R}^{3}}\frac{u_{t}^{2}}{|y|^{\alpha}}dx=t^{1+\alpha}\int_{\mathbb{R}^{3}}\frac{u^{2}}{|y|^{\alpha}}dx.
\end{equation}

\begin{lemma}\label{L:2.3} If $\alpha\in[0,8)$ and $\max\{2,\frac{\alpha+2}{2}\}<p<5$, then there exist $\rho>0$,
$\delta>0$, $e\in H$ with $e\geqslant0$ and $\|e\|_{H}>\rho$ such
that
\begin{description}
 \item $\rm (i)$  $I(u)\geqslant\delta$, for all $u\in H$ with
 $\|u\|_{H}=\rho$.
\item $\rm (ii)$  $I(e)<I(0)=0$.
\end{description}
\end{lemma}
{\it Proof.} (i) Since $H\hookrightarrow L^p(\mathbb{R}^3)$ for
$2\leqslant p<6$, this conclusion is a straightforward consequence
of the definition of $I$.\\

 (ii) For $t>0$ and $u\in H\setminus\{0\}$, by
 (\ref{eq:2.3.1}), (\ref{eq:2.3.2}) and the definition of $I$,
 we see that
\begin{equation}\label{eq:2.3.3}
I(u_t)=\frac{t^3}{2}\|\nabla u\|^2_{2}+\frac{\lambda
t}{2}\|u\|^2_{2}+\frac{t^{1+\alpha}}{2}\int_{\mathbb{R}^{3}}\frac{u^2}{|y|^{\alpha}}dx
+\frac{t^3}{4}\int_{\mathbb{R}^{3}}\phi_u(x)u^2dx
-\frac{t^{2p-1}}{p+1}\|u\|_{p+1}^{p+1}.
\end{equation}
Since $\alpha\in[0,8)$, $p>\max\{2,\frac{\alpha+2}{2}\}$, we see
$I(u_t)\rightarrow-\infty$ as $t\rightarrow+\infty$. Hence, for each
$u\in H\setminus\{0\}$, there is a $t_\ast>0$ large enough such that
(ii) holds with $e=u_{t_\ast}$. Moreover, we may assume that
$e\geqslant0$, otherwise, just replace $e$ by
$|e|$. $\Box$\\
For each $\lambda>0$ and $e$ given by Lemma \ref{L:2.3}, define
\begin{equation}\label{eq:2.4}
c:=c_\lambda=\inf\limits_{\gamma\in\Gamma}\max\limits_{u\in\gamma([0,1])}I_\lambda(u),
\end{equation}
where $\Gamma:=\{\gamma\in C([0,1];H):\gamma(0)=0,\gamma(1)=e\}$.
Clearly, $c>0$ by lemma \ref{L:2.3}. Let $\{t_n\}\subset(0,+\infty)$
be a sequence such that $t_n\rightarrow 1$ as $n\rightarrow
+\infty$, then by (\ref{eq:2.3.1}) it is easy to show that
\begin{equation}\label{eq:2.5}
e_{t_n}:=t_n^2e(t_nx)\rightarrow e \ \text{ in } H,\ \text{as}\
n\rightarrow +\infty.
\end{equation}
Since $I\in C^1(H)$, it follows from Lemma \ref{L:2.3} (ii) that
there is $\varepsilon>0$ small enough such that $I(u)<0$ for all
$u\in B_{\varepsilon}(e)$. Again using (\ref{eq:2.5}), there exists
$t_0\in(0,1)$ such that
\begin{equation}\label{eq:2.5.1}
e_t:=t^2e(tx)\in B_{\varepsilon}(e) \ \text{for all}\ t\in(t_0,1).
\end{equation}
For this $t_0\in(0,1)$, similar to \cite{MBadVBenSRol-reprent} we
have
\begin{lemma}\label{L:2.4} Let $t_0$ be given by (\ref{eq:2.5.1}). Then
for all $t\in(t_0,1)$, we have
\begin{equation}\nonumber
c=\inf\limits_{\gamma\in\Gamma}\max\limits_{u\in\gamma([0,1])}I(u_t)
\end{equation}
where $c$ and $\Gamma$ are defined in (\ref{eq:2.4}),
$u_t=t^2u(tx)$.
\end{lemma}
{\it Proof.} The proof is the same as that of Lemma 11 in
\cite{MBadVBenSRol-reprent}. $\Box$\\

{ By Lemma \ref{L:2.4}}, we know that for any $s\in(t_0,1)$ there
exists $\gamma_s\in\Gamma$ such that
\begin{equation}\label{eq:2.6}
\max\limits_{u\in\gamma_s([0,1])}I(u_s)\leqslant c+(1-s^3).
\end{equation}
For $s\in(t_0,1)$, we define a set
\begin{equation}\label{eq:2.6.1}
U_{s}:=\{u\in\gamma_s([0,1]):I(u)\geqslant c-(1-s^3)\},
\end{equation}
then, (\ref{eq:2.4}) and the definition of $U_{s}$ imply that
$U_{s}\neq\emptyset$ for $s\in(t_0,1)$.

\begin{lemma}\label{L:2.5} If $\alpha\in[0,8)$ and $\max\{2,\frac{\alpha+2}{2}\}<p<5$, then
for $t_0$ given by (\ref{eq:2.5.1}) there exist $t^\ast\in(t_0,1)$
and
$M=\frac{2(c+2)(2p-1)}{(p-2){t^\ast}^3}+\frac{4(c+2)(2p-1)}{(2p-2-\alpha){t^\ast}^{1+\alpha}}$
such that
$$\|u\|^2_{H}+\int_{\mathbb{R}^{3}}\phi_{u}u^2
dx<M\ \ \text{for all}\ u\in U_{s} \text{ with } s\in(t^\ast,1).$$
\end{lemma}
{\bf Proof:} Let $u\in U_{s}$ and note that
$u(x)=(u_s)_\frac{1}{s}$, it follows from (\ref{eq:2.3.1}),
(\ref{eq:2.3.2}) and the definition (\ref{eq:1.10}) that
\begin{eqnarray}\label{eq:2.6.2}
 I(u_s)-I(u)&=\frac{1}{2}(1-\frac{1}{s^3})\|
\nabla u_s\|_{2}^{2}+\frac{\lambda}{2}(1-\frac{1}{s})\|
u_s\|_{2}^{2}+\frac{1}{2}(1-\frac{1}{s^{1+\alpha}})\int_{\mathbb{R}^{3}}\frac{u_s^2}{|y|^\alpha}dx\nonumber\\
&+\frac{1}{4}(1-\frac{1}{s^{3}})\int_{\mathbb{R}^{3}}\phi_{u_s}u_s^2
dx-\frac{1}{p+1}(1-\frac{1}{s^{2p-1}})\| u_s\|_{p+1}^{p+1}.
\end{eqnarray}
For $u\in U_{s}$, (\ref{eq:2.6}) and (\ref{eq:2.6.1}) implies that
\begin{equation}\label{eq:2.6.3}
I(u_s)-I(u)\leqslant2(1-s^3),\ \ \text{for}\ \ s\in(t_0,1).
\end{equation}
By calculation, this and (\ref{eq:2.6.2}) show that, for any $u\in
U_{s}$,
\begin{eqnarray}\nonumber
\frac{\lambda}{2}\frac{s^2-s^3}{s^3-1}\|
u_s\|_{2}^{2}&+&\frac{1}{2}\frac{s^2-s^{3+\alpha}}{s^{3+\alpha}-s^\alpha}\int_{\mathbb{R}^{3}}\frac{u_s^2}{|y|^\alpha}
dx+\frac{1}{p+1}\frac{s^{2p+2}-s^3}{s^{2p+2}-s^{2p-1}}\|
u_s\|_{p+1}^{p+1}\nonumber\\
&-&\frac{1}{2}\| \nabla
u_s\|_{2}^{2}-\frac{1}{4}\int_{\mathbb{R}^{3}}\phi_{u_s}u_s^2
dx\leqslant2s^3.\label{eq:2.7}
\end{eqnarray}
To simplify (\ref{eq:2.7}), we need to use the following facts:
\begin{equation}\nonumber
\frac{s^2-s^{3}}{s^3-1}=\frac{s^2}{s^2+s+1}\geqslant-1\ \
\text{for}\ \ s\geqslant0.
\end{equation}
\begin{equation}\nonumber
g(s)\triangleq\frac{s^2-s^{3+\alpha}}{s^{3+\alpha}-s^\alpha}=\frac{s^{2-\alpha}-s^3}{s^3-1}\overset{s\rightarrow1^-}{\longrightarrow}-\frac{1+\alpha}{3},\
 \text{ and }\ g(s)\equiv g(1)=-1 \text{ if }\ \alpha=2.
\end{equation}
$p>\frac{\alpha+2}{2}$ implies that $\frac{2p-1}{1+\alpha}>1$ and
$\varepsilon_0:=\frac{2p+\alpha}{2(1+\alpha)}\in(1,\frac{2p-1}{1+\alpha})$.
Hence, there is $\delta_1>0$ small enough  such that
$1-\delta_1>t_0$ and
\begin{equation}\nonumber
g(s)\geqslant-\frac{\varepsilon_0(1+\alpha)}{3}=\frac{2p+\alpha}{6}\
\text{for all}\ s\in(1-\delta_1,1).
\end{equation}
Let
\begin{equation}\nonumber
h(s)=\frac{s^{2p+2}-s^3}{s^{2p+2}-s^{2p-1}}=\frac{s^3-s^{4-2p}}{s^3-1}\
\text{for}\ s\in(0,1),
\end{equation}
then
\begin{equation*}\nonumber
\lim\limits_{s\rightarrow 1^-}h(s)=\frac{2p-1}{3}\\ \ \ \
\text{and}\ \
 h{'}(s)=\frac{(2p-1)s^{6-2p}-3s^2-(2p-4)s^{3-2p}}{(s^{3}-1)^2},
\end{equation*}
\begin{equation*}\nonumber
 h{'}(s)\overset{s\rightarrow1^-}{\longrightarrow}-\frac{(2p-1)(p-2)}{3}<0\
\text{if}\ p>2.
\end{equation*}
This shows that there is $\delta_2>0$ small enough and
$1-\delta_2>t_0$ such that
\begin{equation*}\nonumber
 h{'}(s)<0\ \text{and}\ \ h(s)\geqslant\lim\limits_{s\rightarrow
 1^-}h(s)=\frac{2p-1}{3} \ \text{for all}\ s\in(1-\delta_2,1)\
 \text{and}\ p>2.
\end{equation*}
For $p=2$, $h(s)\equiv\frac{2p-1}{3}=1$, so we see that
\begin{equation*}\nonumber
h(s)\geqslant\frac{2p-1}{3} \ \text{for all}\ s\in(1-\delta_2,1)\
 \text{and}\ p\geqslant2.
\end{equation*}
So, for $s\in(t^\ast,1)$ with $t^\ast=1-\min\{\delta_1,\delta_2\}$,
it follows from (\ref{eq:2.7}) that
\begin{eqnarray*}
-\frac{\lambda}{2}\|
u_s\|_{2}^{2}&-&\frac{2p+\alpha}{12}\int_{\mathbb{R}^{3}}\frac{u_s^2}{|y|^\alpha}dx+\frac{1}{p+1}\frac{2p-1}{3}\|
u_s\|_{p+1}^{p+1}\nonumber\\
&-&\frac{1}{2}\| \nabla
u_s\|_{2}^{2}-\frac{1}{4}\int_{\mathbb{R}^{3}}\phi_{u_s}u_s^2
dx\leqslant2s^3.
\end{eqnarray*}
That is
\begin{eqnarray}\label{eq:2.7.1}
-\frac{1}{p+1}\|u_s\|_{p+1}^{p+1}&\geqslant-\frac{3}{2p-1}\left(\frac{\lambda}{2}\|
u_s\|_{2}^{2}+\frac{1}{2}\| \nabla
u_s\|_{2}^{2}+\frac{1}{4}\int_{\mathbb{R}^{3}}\phi_{u_s}u_s^2
dx\right)\nonumber\\
&-\frac{2p+\alpha}{4(2p-1)}\int_{\mathbb{R}^{3}}\frac{u_s^2}{|y|^\alpha}dx-\frac{6}{2p-1}s^3.
\end{eqnarray}
For $u\in U_{s}$, by (\ref{eq:2.6}) it gives that
\begin{equation}\label{eq:2.8}
\begin{split}
&\left(\frac{\lambda}{2}\| u_s\|_{2}^{2}+\frac{1}{2}\| \nabla
u_s\|_{2}^{2}+\frac{1}{4}\int_{\mathbb{R}^{3}}\phi_{u_s}u_s^2
dx\right)\\
&+\frac{1}{2}\int_{\mathbb{R}^{3}}\frac{u_s^2}{|y|^\alpha}dx-\frac{1}{p+1}\|u_s\|_{p+1}^{p+1}\leqslant
c+(1-s^3).
\end{split}
\end{equation}
Hence, it follows from (\ref{eq:2.7.1}) and (\ref{eq:2.8}) that
\begin{equation*}
\begin{split}
\frac{2p-2-\alpha}{4(2p-1)}\int_{\mathbb{R}^{3}}\frac{u_s^2}{|y|^\alpha}dx
&+\frac{2p-4}{2p-1}\left(\frac{\lambda}{2}\|
u_s\|_{2}^{2}+\frac{1}{2}\| \nabla
u_s\|_{2}^{2}+\frac{1}{4}\int_{\mathbb{R}^{3}}\phi_{u_s}u_s^2
dx\right)\\
&\leqslant c+1-\frac{2p-7}{2p-1}s^3\\
&\leqslant c+1+\left|\frac{2p-7}{2p-1}\right|\leqslant c+2\ \text{if
} p>2\ \text{and}\ s<1.
\end{split}
\end{equation*}
This implies that, if $5>p>\max\{2,\frac{\alpha+2}{2}\}$ and
$s\in(t^\ast,1)$
\begin{equation}\label{eq:2.9}
\frac{\lambda}{2}\| u_s\|_{2}^{2}+\frac{1}{2}\| \nabla
u_s\|_{2}^{2}+\frac{1}{4}\int_{\mathbb{R}^{3}}\phi_{u_s}u_s^2
dx\leqslant \frac{(c+2)(2p-1)}{2(p-2)}.
\end{equation}
and
\begin{equation*}
\frac{1}{4}\int_{\mathbb{R}^{3}}\frac{u_s^2}{|y|^\alpha}dx\leqslant
\frac{(c+2)(2p-1)}{2p-2-\alpha}\ \ \text{for}\ \alpha\in[0,8).
\end{equation*}
Hence, it follows from (\ref{eq:2.9}) and by using (\ref{eq:2.3.1})
and (\ref{eq:2.3.2})
\begin{equation*}
\frac{\lambda}{2}s\| u\|_{2}^{2}+\frac{1}{2}s^3\| \nabla
u\|_{2}^{2}+\frac{1}{4}s^3\int_{\mathbb{R}^{3}}\phi_{u}u^2
dx\leqslant \frac{(c+2)(2p-1)}{2(p-2)}.
\end{equation*}
Since $s\in (t^\ast,1)$, $s\geq s^3\geqslant {t^\ast}^3$ and
$s^{1+\alpha}\geqslant {t^\ast}^{1+\alpha}$ for $\alpha\in[0,8)$,
those and $p>\max\{2,\frac{\alpha+2}{2}\}$ imply that
\begin{equation}\label{eq:2.9.1}
\|u\|^2_{H}+\int_{\mathbb{R}^{3}}\phi_{u}u^2
dx\leqslant\frac{2(c+2)(2p-1)}{(p-2){t^\ast}^3}+\frac{4(c+2)(2p-1)}{(2p-2-\alpha){t^\ast}^{1+\alpha}},\
\
\end{equation}
and Lemma \ref{L:2.5} is proved by taking
$M=\frac{2(c+2)(2p-1)}{(p-2){t^\ast}^3}+\frac{4(c+2)(2p-1)}{(2p-2-\alpha){t^\ast}^{1+\alpha}}$. $\Box$\\

Note that $M$ given by the above lemma depends on $\lambda$, since
$c$ depends on $\lambda$ by the definition of $I$. The following
lemma is for getting a bounded (PS) sequence. In this lemma, the
constant $M$ can be chosen independent of $\lambda$ if
$\lambda\in(0,1]$.
\begin{lemma}\label{L:2.6} If $\alpha\in[0,8)$, $\max\{2,\frac{\alpha+2}{2}\}<p<5$ and $c$ be given by (\ref{eq:2.4}). Then
there exists a bounded nonnegative sequence $\{u_n\}\subset H$ such
that
\begin{equation}\label{eq:2.10}
I(u_n)\rightarrow c>0,\ \ I{'}(u_n)\rightarrow0 \ \text{as}\
n\rightarrow+\infty,
\end{equation}
Moreover, if $\lambda\in(0,1]$ there exists $M>0$ which is
independent of $\lambda\in(0,1]$ such that
$$\|u_n\|^2_{H}+\int_{\mathbb{R}^{3}}\phi_{u_n}{u_n}^2
dx\leqslant M.$$
\end{lemma}
{\bf Proof:} For $t\in(t^\ast,1)$ with $t^\ast$ given in Lemma
\ref{L:2.5}, let
\begin{equation}\label{eq:2.10.1}
W_{t}=\{|u|:u\in U_{t}\}, \ \ U_{t} \text{  defined in
(\ref{eq:2.6.1}),}
\end{equation}
and then for $u\in W_{t}$, by (\ref{eq:2.9.1}) (\ref{eq:2.3.3}) and
(\ref{eq:2.6}) we have that
\begin{equation*}
\begin{split}
 I(u)-I(u_t)&=\frac{1}{2}(1-t^3)\|\nabla
u\|_{2}^{2}+\frac{\lambda}{2}(1-t)\|
u\|_{2}^{2}+\frac{1-t^{1+\alpha}}{2}\int_{\mathbb{R}^{3}}\frac{u^2}{|y|^\alpha}dx\nonumber\\
&+\frac{1-t^3}{4}\int_{\mathbb{R}^{3}}\phi_{u}u^2
dx-\frac{1-t^{2p-1}}{p+1}\| u\|_{p+1}^{p+1}\\
&\leqslant(1-t^3)\left(\frac{\lambda}{2}\| u\|_{2}^{2}+\frac{1}{2}\|
\nabla u\|_{2}^{2}+\frac{1}{4}\int_{\mathbb{R}^{3}}\phi_{u}u^2
dx\right)\\
&+\frac{1-t^{2p-1}}{t^{2p-1}}\left(\frac{1}{2}\int_{\mathbb{R}^{3}}\frac{u_t^2}{|y|^\alpha}dx-\frac{1}{p+1}\| u_t\|_{p+1}^{p+1}\right)\\
&\leqslant(1-t^3)\left(\frac{\lambda}{2}\| u\|_{2}^{2}+\frac{1}{2}\|
\nabla u\|_{2}^{2}+\frac{1}{4}\int_{\mathbb{R}^{3}}\phi_{u}u^2
dx\right)+\frac{1-t^{2p-1}}{t^{2p-1}}I(u_t)\\
&\leqslant(1-t^3)\frac{(c+2)(2p-1)}{(2p-4){t^\ast}^3}+(1-t^{2p-1})\frac{c+1}{{t^\ast}^{2p-1}}
\rightarrow 0\ \text{as}\ t\rightarrow1^-.
\end{split}
\end{equation*}
 On the other hand, similar to (\ref{eq:2.6.3}) we know that
\begin{equation*}
\begin{split}
 I(u_t)-I(u)\leqslant2(1-t^3)
\rightarrow 0\ \text{as}\ t\rightarrow1^-.
\end{split}
\end{equation*}
Hence,
\begin{equation}\label{eq:2.11}
 \limsup\limits_{t\rightarrow1^-\ u\in W_{t}}|I(u_t)-I(u)|=0.
\end{equation}
Define
\begin{equation*}
 S=\left\{|u|:u\in H\ \text{and}\ \|u\|^2_{H}+\int_{\mathbb{R}^{3}}\phi_{u}{u}^2
dx<M\right\},\text{ where} \ M \text{ is given by Lemma \ref{L:2.5}},
\end{equation*}
\begin{equation}\label{eq:2.11.1}
 S_\delta=\left\{u:u\in H\ \text{and}\
 dist(u,S)<\delta\right\},\ \delta\in(0,1).
\end{equation}
Clearly, $\|\upsilon\|_{H}\leqslant \sqrt{M}+1$ for all $\upsilon\in
S_\delta$. Then, by Lemma \ref{L:2.2}, there is a constant $K:=K(M)$
such that
\begin{equation}\label{eq:2.12}
 \|I{'}(u)-I{'}(\upsilon)\|_{H'}\leqslant
K\|u-\upsilon\|_{H}\ \text{for all } u,\upsilon\in S_\delta.
\end{equation}
and since $I\in C^{1}(H,\mathbb{R})$, there exists $C_{S}>0$ such
that
\begin{equation}\label{eq:2.12.0}
 \|I(u)-I(\upsilon)\|_{H}\leqslant
C_{S}\|u-\upsilon\|_{H}\ \text{for all } u,\upsilon\in S_\delta.
\end{equation}
For any $\ m\in \mathbb{N}$ and $M$ given by Lemma \ref{L:2.5}, let
\begin{equation}\label{eq:2.12.1}
 \Lambda_{m}=\left\{|u|:u\in H,\
 \|u\|^2_{H}+\int_{\mathbb{R}^{3}}\phi_{u}{u}^2
dx<M+\frac{1}{m}\text{ and
 }|I(u)-c|\leqslant\frac{C_{S}+1}{\sqrt{m}}\right\}.
\end{equation}
We claim that  $\Lambda_{m}\neq\emptyset$.
Indeed, for any $m\geqslant1$, since (\ref{eq:2.11}) we can find
$t_m\in(t^\ast,1)$ such that
\begin{equation*}
1-{t_m}^3<\frac{1}{32m}\text{ and
 }I(u)\leqslant I(u_{t_m})+\frac{1}{32m}\text{ for all }u\in
 W_{t_m}.
\end{equation*}
Then it follows from (\ref{eq:2.6}) and (\ref{eq:2.6.1}) that
\begin{equation}\label{eq:2.12.2}
c-\frac{1}{32m}\leqslant I(u)\leqslant c+\frac{1}{16m}\text{ for all
}u\in
 W_{t_m}.
\end{equation}
By the definition of $W_{t_m}$, Lemma \ref{L:2.5} implies that
$\|u\|^2_{H}+\int_{\mathbb{R}^{3}}\phi_{u}{u}^2 dx\leqslant M$ for
all $u\in W_{t_m}$. This and (\ref{eq:2.12.2}) show that
$W_{t_m}\subset\Lambda_{m}$, that is $\Lambda_{m}\neq\emptyset$.

Next, we claim that there are infinitely many elements in
$\{\Lambda_{m}\}_{m=1}^{+\infty}$, which we still simply denote by
$\Lambda_{m}$ ($m=1,2,\cdots$,), such that for each $m\geqslant1$,
there is $u_m\in\Lambda_{m}$ with
\begin{equation}\label{eq:2.13}
 \|I{'}(u_m)\|_{H'}< \frac{1+K}{\sqrt{m}},\ K\text{ is
given by (\ref{eq:2.12})}.
\end{equation}
Then, to prove Lemma \ref{L:2.6}
 we need only to show the above claim. By contradiction, if the
 claim is false, then there must be a number $\bar{m}\in\mathbb{N}$
 with $\bar{m}>\max\{\frac{1}{8c},4\}$ such that
\begin{equation}\label{eq:2.13.1}
 \|I{'}(u)\|_{H'}\geqslant
\frac{1+K}{\sqrt{m}},\text{ for all }m>\bar{m}\text{ and }
u\in\Lambda_{m} .
\end{equation}
By the above discussion we know that $W_{t_m}\subset\Lambda_{m}$.
For any $u\in W_{t_m}$, the definition of $W_{t_m}$ and Lemma
\ref{L:2.5} show that
$\|u\|^2_{H}+\int_{\mathbb{R}^{3}}\phi_{u}{u}^2 dx\leqslant M$ and
$W_{t_m}\subset S$. Hence,
\begin{equation}\nonumber
 W_{t_m}\subset S\cap\{u\in
H:|I(u)-c|<\frac{1}{8m}\}\subset S\cap\{u\in
H:|I(u)-c|<\frac{C_{S}+1}{\sqrt{m}}\}\subset \Lambda_{m},
\end{equation}
 where
(\ref{eq:2.12.2}) is used. Then $$S\cap\{u\in
H:|I(u)-c|<\frac{C_{S}+1}{\sqrt{m}}\}\neq\emptyset.$$ Let
$\varepsilon=\frac{1}{16m}$,  $\delta=\frac{1}{2\sqrt{m}}$, then
$\frac{8\varepsilon}{\delta}=\frac{1}{\sqrt{m}}<\frac{1}{2}<1$,
since $\bar{m}>\max\{\frac{1}{8c},4\}$. So,
\begin{equation*}
 (S)_{2\delta}=S_{\frac{1}{\sqrt{m}}}=\left\{u:u\in H\ \text{and}\
 dist(u,S)<\frac{1}{\sqrt{m}}\right\}.
\end{equation*}
By the definitions of $S$ and $\Lambda_m$, we have
$$S\cap\{u\in H:|I(u)-c|<\frac{C_{S}+1}{\sqrt{m}}\}\subset
\Lambda_{m}.$$
Hence, for any $u\in S\cap\{u\in
H:|I(u)-c|<\frac{C_{S}+1}{\sqrt{m}}\}\subset \Lambda_{m}$,
\begin{equation}\label{eq:2.14}
 \|I{'}(u)\|_{H'}\geqslant
\frac{1+K}{\sqrt{m}},\text{ for all }m>\bar{m}.
\end{equation}
For any $\upsilon\in
 S_{\frac{1}{\sqrt{m}}}\cap\{u\in H:|I(u)-c|<\frac{1}{8m}\}$, it is not difficult to
know that there is $u_0\in S$ such that
\begin{equation}\label{eq:2.14.1}
 \|u_0-\upsilon\|_{H}<\frac{1}{\sqrt{m}}.
\end{equation}
This and (\ref{eq:2.12.0}) show that
\begin{eqnarray}\nonumber
\|I(u_0)-c\|_{H}&\leqslant&\|I(\upsilon)-I(u_0)\|_{H}+\|I(\upsilon)-c\|_{H}\nonumber\\
&\leqslant&\|I(\upsilon)-c\|_{H}+\frac{C_{S}}{\sqrt{m}}\nonumber\\
&\leqslant&\frac{1}{8m}+\frac{C_{S}}{\sqrt{m}}\leqslant\frac{C_{S}+1}{\sqrt{m}}\nonumber.
\end{eqnarray}
That is $u_0\in S\cap\{u\in H:|I(u)-c|<\frac{C_{S}+1}{\sqrt{m}}\}$.
Then, it follows from (\ref{eq:2.12}), (\ref{eq:2.14}) and
(\ref{eq:2.14.1}) that, for $\upsilon\in
S_{\frac{1}{\sqrt{m}}}\cap\{u\in H:|I(u)-c|<\frac{1}{8m}\}$,
\begin{eqnarray}\nonumber
\|I{'}(\upsilon)\|_{H'}&=&\|I{'}(\upsilon)-I{'}(u_0)+I{'}(u_0)\|_{H'}\nonumber\\
&\geqslant&\|I{'}(u_0)\|_{H'}-\|I{'}(\upsilon)-I{'}(u_0)\|_{H'}\nonumber\\
&\geqslant&\frac{1+K}{\sqrt{m}}-K\|u_0-\upsilon\|_{H}\nonumber\\
&\geqslant&\frac{1+K}{\sqrt{m}}-K\frac{1}{\sqrt{m}}=\frac{1}{\sqrt{m}}.\nonumber
\end{eqnarray}
Applying Lemma \ref{L:2.1} with $X=H$, $\varphi=I$, we know that
there is an homeomorphism $\eta(t,\cdot):[0,1]\times H\rightarrow H$
such that
\begin{equation}\label{i}
 \eta(t,u)=u, \text{ if } t=0\text{ or }
 u\notin S_{\frac{1}{\sqrt{m}}}\cap\{u\in H:|I(u)-c|\leqslant\frac{1}{8m}\};
\end{equation}
\begin{equation}\label{ii}
 I(\eta(1,u))\leqslant c-\frac{1}{16m},\text{ for } u\in S\cap\{u\in H:|I(u)-c|\leqslant\frac{1}{8m}\};
\end{equation}
\begin{equation}\label{iii}
I(\eta(t,u))\leqslant I(u), \text{ for any } u\in H.
\end{equation}
Let $\xi(u):=\eta(1,u)$ and
$\bar{\gamma}(t)=\xi(|\gamma_{t_m}(t)|)\in C([0,1],H).$ By
$m>\bar{m}>\max\{\frac{1}{8c},4\}$, $c>\frac{1}{8m}$, then
$\{0,e\}\nsubseteq S_{\frac{1}{\sqrt{m}}}\cap\{u\in
H:|I(u)-c|<\frac{1}{8m}\}$, since $I(e)<0$ and $|I(e)-c|=c+|I(e)|>c$
where $e$ is given by Lemma \ref{L:2.3}. With this observation and
(\ref{i}) we see that
$\bar{\gamma}(0)=\xi(|\gamma_{t_m}(0)|)=\xi(0)=\eta(1,0)=0$,
$\bar{\gamma}(1)=\xi(|\gamma_{t_m}(1)|)=\xi(e)=\eta(1,e)=e$. Hence,
 $\bar{\gamma}\in\Gamma$, with $\Gamma$ defined in
(\ref{eq:2.4}). For each $m\geqslant \bar{m}$, let
$u_m\in\bar{\gamma}([0,1])$ be such that
\begin{equation}\label{eq:2.15}
 I(\xi(|u_m|)=\max\limits_{u\in\gamma_{t_m}[0,1]}I(\xi(|u|))=\max\limits_{v\in\bar{\gamma}[0,1]}I(v)
\geqslant c.
\end{equation}
Since $u_m\in\gamma_{t_m}[0,1]$,
$|u_m|\in|\gamma_{t_m}[0,1]|=\{|u|:u\in\gamma_{t_m}[0,1]\}$. We are
ready to get a contradiction in both of the following two
cases.\\
{\bf Case A:} If $|u_m|\in|\gamma_{t_m}[0,1]|\setminus U_{t_m}$,
then (\ref{iii}) and the definition of $U_{t_m}$ imply that
\begin{equation*}
 I(\xi(|u_m|)=I(\eta(1,|u_m|))\leqslant I(u_m)\leqslant
 c-(1-t_m^3)<c,
\end{equation*}
which contradicts (\ref{eq:2.15}).\\
{\bf Case B:} If $|u_m|\in U_{t_m}$, then by (\ref{eq:2.10.1})
$|u_m|\in W_{t_m}$ and (\ref{eq:2.12.2}) implies that
$|I(|u_m|)-c|\leqslant \frac{1}{16m}$. Moreover,
$\|u_m\|^2_{H}+\int_{\mathbb{R}^{3}}\phi_{u_m}{u_m}^2 dx\leqslant M$
by Lemma \ref{L:2.5}. Hence $|u_m|\in S\cap\{u\in
H:|I(u)-c|\leqslant\frac{1}{16m}\}$, and it follows from (\ref{ii})
that
\begin{equation*}
 I(\xi(|u_m|)=I(\eta(1,|u_m|))\leqslant
 c-\frac{1}{16m}< c,
\end{equation*}
this is a contradiction to (\ref{eq:2.15}). $\Box$\\

\section{Existence for $\lambda>0$:  Proof of Theorem \ref{th1.1}.}
Motivated by \cite{BadMBenciVRolS-JEMS}, we prove Theorem
\ref{th1.1} by a result due to S.Solimini \cite{SSol-AnnIHPA}, which
is a version of so called concentration-compactness principle. To
state this result, we should recall the operator $T_{s,\xi}$ and its
basic properties. Let $s>0$, $N\geq3$ and $\xi\in \mathbb{R}^N$ be
fixed, for any $u\in L^q(\mathbb{R}^N)$ ($1<q<+\infty$) we define
\begin{equation}\label{eq:3.1}
 T_{s,\xi}u(x)\triangleq
 T(s,\xi)u(x):=s^{-\frac{N-2}{2}}u(s^{-1}x+\xi),\
 \forall x\in\mathbb{R}^N.
\end{equation}
Clearly, $T(s,\xi)u\in L^{q}(\mathbb{R}^N)$ if $u\in
L^{q}(\mathbb{R}^N)$ and $T(s,\xi)$ is also well defined on Hilbert
space $D^{1,2}(\mathbb{R}^N)$ with scalar product
\begin{equation}\label{insc}
\langle u,v\rangle=\int_{\mathbb{R}^N}\nabla u\nabla vdx, \text{ for
} u,v\in D^{1,2}(\mathbb{R}^N),
\end{equation}
 since $T(s,\xi)u\in
D^{1,2}(\mathbb{R}^N)$ if $u\in D^{1,2}(\mathbb{R}^N)$. It is not
difficult to see that the linear operators
\begin{equation*}
 u\in L^{2^\ast}(\mathbb{R}^N)\rightarrowtail T(s,\xi)u\in L^{2^\ast}(\mathbb{R}^N)
 \ \text{and } u\in D^{1,2}(\mathbb{R}^N)\rightarrowtail T(s,\xi)u\in D^{1,2}(\mathbb{R}^N)
\end{equation*}
are isometric, where $2^\ast=\frac{2N}{N-2}$. Moreover, we have that
\begin{equation}\label{eq:3.2}
T_{s,\xi}^{-1}=T(s^{-1},-s\xi), \quad\quad
T_{s,\xi}T_{\mu,\eta}=T(s\mu,\xi/\mu+\eta).
\end{equation}
\begin{equation}\label{eq:3.3}
\|\nabla T_{s,\xi}u\|_{2}^{2}=\|\nabla u\|_{2}^{2},\ \ \
\|T_{s,\xi}u\|_{q}^{q}={s}^{N-\frac{q(N-2)}{2}}\|u\|_{q}^{q}.
\end{equation}
For $N\geqslant3$, $k\in[2,N)$ and $x\in \mathbb{R}^N$, in this
section we denote that
$$x=(y,z)\in\mathbb{R}^k\times\mathbb{R}^{N-k}, \text{ i.e. }
y\in\mathbb{R}^k, z\in\mathbb{R}^{N-k},$$
$\tilde{y}=(y,0)\in\mathbb{R}^k\times\mathbb{R}^{N-k}$,
$\tilde{z}=(0,z)\in\mathbb{R}^k\times\mathbb{R}^{N-k}$. Similarly,
$x_n=(y_n,z_n)\in\mathbb{R}^k\times\mathbb{R}^{N-k}$,
$\tilde{y}_n=(y_n,0)\in\mathbb{R}^k\times\mathbb{R}^{N-k}$.

\begin{lemma}\label{L:3.1}
             (\cite{BadMBenciVRolS-JEMS}, Proposition 22) Let $\{\eta_{n}\}\subset\mathbb{R}^{N}$ be such that
$\lim\limits_{n\rightarrow\infty}|\eta_{n}|=\infty$ and fix $R>0$.
Then for any $m\in\mathbb{N}\setminus\{0,1\}$ there exists
$N_m\in\mathbb{N}$ such that for any $n>N_m$ one can find a sequence
of unit orthogonal matrices, $\{g_i\}_{i=1}^{m}\in O(N)$ satisfying
the condition
\begin{eqnarray}\nonumber
B_{R}(g_i\eta_n)\cap B_{R}(g_j\eta_n)=\emptyset,
\quad\text{for}\quad i\neq j.
\end{eqnarray}
\end{lemma}
\begin{lemma}\label{L:3.2}
(\cite{BadMBenciVRolS-JEMS}, Proposition 11) Let $q\in(1,\infty)$
and $\{s_n\}\subset(0,\infty)$, $\{\xi_n\}\subset \mathbb{R}^N$ be
such that $s_n\overset{n}{\rightarrow}s\neq0$,
$\xi_n\overset{n}{\rightarrow}\xi$. Then
\begin{eqnarray}\nonumber
T_{s_n,\xi_n}u_n\overset{n}{\rightharpoonup}T_{s,\xi}u \text{ weakly
in } L^{q}(\mathbb{R}^N),
\end{eqnarray}
if $u_n\overset{n}{\rightharpoonup}u$ weakly in
$L^{q}(\mathbb{R}^N)$.
\end{lemma}
\begin{lemma}\label{L:3.2.0}
 Let $\{s_n\}\subset(0,\infty)$, $\{\xi_n\}\subset
\mathbb{R}^N$ be such that $s_n\overset{n}{\rightarrow}s_0\neq0$,
$\xi_n\overset{n}{\rightarrow}\xi$. If
$v_n\overset{n}{\rightharpoonup}v$ weakly in
$D^{1,2}(\mathbb{R}^N)$, then
\begin{eqnarray}\nonumber
T_{s_n,\xi_n}v_n\overset{n}{\rightharpoonup}T_{s_0,\xi}v \text{
weakly in } D^{1,2}(\mathbb{R}^N).
\end{eqnarray}
\end{lemma}
{\it Proof.} For any $\varphi\in C^{\infty}_{0}(\mathbb{R}^{N})$, by
(\ref{insc}) we get that
\begin{equation}\label{eq:3.2.0.1}
 \langle T^{-1}_{s_n,0}v_n,\varphi\rangle=\langle
 v_n,T_{s_n,0}\varphi\rangle=\langle
 v_n,T_{s_0,0}\varphi\rangle+\langle
 v_n,T_{s_n,0}\varphi-T_{s_0,0}\varphi\rangle.
\end{equation}
Since
\begin{eqnarray*}
\lim\limits_{n\rightarrow\infty}
\|\nabla(T_{s_n,0}\varphi-T_{s_0,0}\varphi)\|^2_{2}
&=&\lim\limits_{n\rightarrow\infty}\int_{\mathbb{R}^{N}}|\nabla
T_{s_n,0}\varphi|^2dx+\int_{\mathbb{R}^{N}}|\nabla
T_{s_0,0}\varphi|^2dx\\&-&2\lim\limits_{n\rightarrow\infty}\int_{\mathbb{R}^{N}}\nabla
T_{s_n,0}\varphi\nabla T_{s_0,0}\varphi=0,\\
\end{eqnarray*}
we have
\begin{equation}\label{eq:3.2.0.2}
 \langle v_n,T_{s_n,0}\varphi-T_{s_0,0}\varphi\rangle
 \leqslant\|\nabla v_n\|_{2}\|\nabla(T_{s_n,0}\varphi-T_{s_0,0}\varphi)\|_{2}\overset{n}{\rightarrow}0.
\end{equation}
By $T_{s_0,0}\varphi\in C_0^\infty(\mathbb{R}^N)$ and
$v_n\overset{n}{\rightharpoonup}v$ weakly in
$D^{1,2}(\mathbb{R}^N)$, we have
\begin{equation}\label{eq:3.2.0.3}
\langle
 v_n,T_{s_0,0}\varphi\rangle\overset{n}{\longrightarrow}\langle
 v,T_{s_0,0}\varphi\rangle=\langle T^{-1}_{s_0,0}
 v,\varphi\rangle.
\end{equation}
It follows from (\ref{eq:3.2.0.1}) to (\ref{eq:3.2.0.3}) that
\begin{equation}\label{eq:3.2.0.4}
 \langle T^{-1}_{s_n,0}v_n,\varphi\rangle\overset{n}{\longrightarrow}\langle
 T^{-1}_{s_0,0}
 v,\varphi\rangle ,\text{ for any $\varphi\in C_0^\infty(\mathbb{R}^N)$}.
\end{equation}
On the other hand,  for any $\psi\in D^{1,2}(\mathbb{R}^{N})$ and
any $\epsilon>0$, there exists $\varphi\in
C^{\infty}_{0}(\mathbb{R}^{N})$ such that
$\|\nabla(\psi-\varphi)\|_{2}<\epsilon$ and
\begin{equation}\nonumber
 \langle
 T^{-1}_{s_n,0}v_n,\psi-\varphi\rangle\leqslant\|\nabla(T^{-1}_{s_n,0}v_n)\|_{2}\|\nabla(\psi-\varphi)\|_{2}
 =\|\nabla v_n\|_{2}\|\nabla(\psi-\varphi)\|_{2},
\end{equation}
 this and (\ref{eq:3.2.0.4}) imply that
\begin{equation}\nonumber
 \langle T^{-1}_{s_n,0}v_n,\varphi\rangle\overset{n}{\longrightarrow}\langle T^{-1}_{s_0,0}
 v,\varphi\rangle,\text{ for any } \varphi\in
 D^{1,2}(\mathbb{R}^{N}).\ \ \ \ \Box
\end{equation}
\begin{lemma}\label{L:3.3}
(\cite{SSol-AnnIHPA}, A corollary of Theorem 1) If $\{u_n\}\subset
D^{1,2}(\mathbb{R}^N)$ is bounded, then, up to a subsequence, either
$u_n\overset{n}{\rightarrow}0$ in $L^{2^\ast}(\mathbb{R}^N)$ or
there exist  $\{s_n\}\subset(0,\infty)$ and $\{\xi_n\}\subset
\mathbb{R}^N$ such that
\begin{eqnarray}\nonumber
T_{s_n,\xi_n}u_n\overset{n}{\rightharpoonup}u\neq0 \text{ weakly in
} L^{2^\ast}(\mathbb{R}^N).
\end{eqnarray}
\end{lemma}
Let
$$D_s^{1,2}(\mathbb{R}^{N})\overset{\triangle}{=}\{u\in
D^{1,2}(\mathbb{R}^{N}): u(x)=u(y,z)=u(|y|,z)\},$$ we see that
$D_s^{1,2}(\mathbb{R}^{N})\subset D^{1,2}(\mathbb{R}^{N})$ is a
closed set, hence $D_s^{1,2}(\mathbb{R}^{N})$ is a Hilbert space
with scalar product as (\ref{insc}). Based on Lemmas \ref{L:3.1} to
\ref{L:3.3}, we have the following lemma which ensures us to get a
nontrivial solution for (\ref{eq:1.1}) without proving the (PS)
condition.
\begin{lemma}\label{L:3.4}
If $\{u_n\}\subset D_s^{1,2}(\mathbb{R}^{N})$ is bounded and there
exist $\{s_{n}\}\subset(0,+\infty)$ and
$\{x_{n}\}\subset\mathbb{R}^{N}$ with
$x_n=(y_n,z_n)\in\mathbb{R}^{k}\times\mathbb{R}^{N-k}$ such that
\begin{eqnarray}\label{eq:3.4.a}
T(s_{n},x_{n})u_{n}\overset{n}{\rightharpoonup} u\neq0 \text{ weakly
in } L^{2^{\ast}}(\mathbb{R}^{N}).
\end{eqnarray}
 Then
 \begin{eqnarray}\nonumber
v_n=T(s_n,0)w_{n}\overset{n}{\rightharpoonup} v\not\equiv0 \text{
weakly in } D_s^{1,2}(\mathbb{R}^{N}),
\end{eqnarray}
where $w_n=T(1,\tilde{z}_{n})u_{n}$ and $\tilde{z}_n=(0,z_n)$.
Moreover, if $\{u_n\}$ is also bounded in $L^q(\mathbb{R}^N)$ for
some $1<q<2^\ast$, then, there exists constant $l>0$ such that
$s_{n}>l$ for all $n$.
\end{lemma}
{\bf Proof:} The proof of this lemma is almost the same as that of
Lemma 23 in \cite{BadMBenciVRolS-JEMS}. But for the sake of
completeness, we give its proof.\\
  Since $\{u_n\}$ is bounded in $D_s^{1,2}(\mathbb{R}^{N})$, by the
  definition of $T_{s,\xi}$ we see that $\{v_n\}$ is also
  bounded in $D_s^{1,2}(\mathbb{R}^{N})$. Then there is $v\in
  D_s^{1,2}(\mathbb{R}^{N})$ such that
 \begin{eqnarray}\nonumber
v_n=T(s_n,0)w_{n}\overset{n}{\rightharpoonup} v \text{ weakly in }
D_s^{1,2}(\mathbb{R}^{N}).
\end{eqnarray}
We claim that $v\not\equiv0$. Otherwise if $v\equiv0$, then it leads
to a contradiction in the following two cases. For $x_n=(y_n,z_n)$,
we note that
\begin{eqnarray}\nonumber
\tilde{y}_n=(y_n,0)\in\mathbb{R}^k\times\mathbb{R}^{N-k},\ \ \
\tilde{z}_n=(0,z_n)\in\mathbb{R}^k\times\mathbb{R}^{N-k}.
\end{eqnarray}
{\bf Case A:} If $\{s_n\tilde{y}_n\}\subset\mathbb{R}^N$ is bounded.
Then, there is
$\tilde{y}_0=(y_0,0)\in\mathbb{R}^k\times\mathbb{R}^{N-k}$ such that
$s_n\tilde{y}_n\overset{n}{\longrightarrow}\tilde{y}_0$ and from
(\ref{eq:3.2}) we have
$$T_{1,-s_{n}\tilde{y}_n}T_{s_{n},\tilde{y}_n}w_{n}=
T_{1,-s_{n}\tilde{y}_n}T_{s_{n},x_n}u_n\overset{n}{\rightharpoonup}
T_{1,-\tilde{y}_0}u\neq0\quad\text{in}\quad
L^{2^{\ast}}(\mathbb{R}^{N}),$$ where we have used assumption
(\ref{eq:3.4.a}) and Lemma \ref{L:3.2}. On the other hand, since
$v\equiv0$, from (\ref{eq:3.2}) we have
$$T_{1,-s_{n}\tilde{y}_n}T_{s_{n},\tilde{y}_n}w_{n}
=T_{s_{n},0}w_n=v_n\overset{n}{\rightharpoonup}0\quad\text{in}\quad
D^{1,2}(\mathbb{R}^{N}),$$ then we have a contradiction.\\
 {\bf Case B:} If $|s_n\tilde{y}_n|\rightarrow+\infty$. We claim that there is also a contradiction. Indeed, since $u\not\equiv0$,
there exist $\Omega\subset\mathbb{R}^{N}$, $|\Omega|\neq0$ and
$\kappa>0$ such that $u>\kappa$ or $u<-\kappa$ a.e in $\Omega$. So
we can choose $R>0$ such that $|B_R\cap\Omega|>0$ and
$$\left|\int_{\mathbb{R}^{N}}T_{s_n,\tilde{y}_n}w_n\chi_{B_R\cap\Omega}dx\right|
\overset{n}{\rightarrow}\left|\int_{\mathbb{R}^{N}}u\chi_{B_R\cap\Omega}dx\right|\geqslant\kappa|B_R\cap\Omega|>0.$$
But,
$$T_{s_n,\tilde{y}_n}w_n=T_{s_n,\tilde{y}_n}T_{s_n^{-1},0}v_n=T_{1,s_n\tilde{y}_n}v_n.$$
Then,
\begin{eqnarray}\nonumber
\left|\int_{\mathbb{R}^{N}}T_{s_n,\tilde{y}_n}w_n\chi_{B_R\cap\Omega}dx\right|
&\leqslant&\int_{B_R}|T_{s_n,\tilde{y}_n}w_n|dx\nonumber\\
&=&\int_{B_R(s_n\tilde{y}_n)}|v_n|dx\nonumber\\
&\leqslant&
C_R\left\{\int_{B_R(s_n\tilde{y}_n)}|v_n|^{2^{\ast}}dx\right\}^{\frac{1}{2^{\ast}}}.\nonumber
\end{eqnarray}
This implies
$$\inf\limits_{n}\int_{B_R(s_n\tilde{y}_n)}|v_n|^{2^{\ast}}dx>\epsilon>0.$$
Since $|s_n\tilde{y}_n|\rightarrow+\infty$, by Lemma \ref{L:3.2} we
have that for any $m\in\mathbb{N}$ we have $\{g_i\}_{i=1}^{m}\subset
O(N)$ and $n_m\in\mathbb{N}$ such that
\begin{eqnarray}\nonumber
\int_{\mathbb{R}^{N}}|u_n|^{2^{\ast}}dx&=&\int_{\mathbb{R}^{N}}|v_n|^{2^{\ast}}dx\geqslant
\sum_{i=1}^{m}\int_{B_R(g_i(s_n\tilde{y}_n))}|v_n|^{2^{\ast}}dx\nonumber\\
&=&m\int_{B_R(s_n\tilde{y}_n)}|v_n|^{2^{\ast}}dx>m\epsilon\quad\text{for}\quad
n>n_m,\nonumber
\end{eqnarray}
 where we have used (\ref{eq:3.3}) and $v(y,z)=v(|y|,z)$. Let
$m\rightarrow\infty$, we have
$\|u_n\|_{2^\ast}\overset{n}{\longrightarrow}+\infty$,
 which contradicts that $\{u_n\}\subset L^{2^\ast}
$ is bounded. \\

Now we can choose $\varphi\in C_{0}^{\infty}(\mathbb{R}^{N})$
satisfing $\int_{\mathbb{R}^{N}}v\varphi dx\neq0$. Choose $R>0$ such
that $supp\varphi\subset B_R$. Since $u\in
D^{1,2}(\mathbb{R}^{N})\rightarrow T(s,\xi)u\in
D^{1,2}(\mathbb{R}^{N})$ is isometric, we obtain
$\{T_{\lambda_n,0}w_n\}$ is bounded in $D^{1,2}(B_R)$, so is in
$L^{2}(B_R)$, hence $T_{s_n,0}w_n\rightharpoonup v$ in $L^{2}(B_R)$.
Then we have
$$\int_{\mathbb{R}^{N}}T_{s_n,0}w_n\varphi dx
=\int_{B_R}T_{s_n,0}w_n\varphi dx\rightarrow\int_{B_R}v\varphi
dx=\int_{\mathbb{R}^{N}}v\varphi dx\neq0$$ On the otherhand we have
\begin{eqnarray}\nonumber
\left|\int_{\mathbb{R}^{N}}T_{s_n,0}w_n\varphi dx\right|
&\leqslant&\|\varphi\|_\infty|B_R|^{\frac{q-1}{q}}\|T_{s_n,0}w_n\|_{L^q(B_R)}\nonumber\\
&\leqslant&s^{\frac{N}{q}-\frac{N-2}{2}}_n\|\varphi\|_{\infty}|B_R|^{\frac{q-1}{q}}\sup\limits_{n}\|u_n\|_{q}.\nonumber
\end{eqnarray}
Since $1<q<2^\ast$, $\frac{N}{q}-\frac{N-2}{2}>0$. So, if
$\lim\limits_{n\rightarrow\infty}s_n=0$, we obtain a contradiction.
This implies that there exists $l>0$ such that
$\inf\limits_{n} s_{n}>l$, since $s_n>0$ for all $n$.$\Box$\\

\begin{lemma}\label{Le:gz}
Let $u\in D^{1,2}(\mathbb{R}^{N})\setminus\{0\}$ be a nonnegative
function, and $K\subset \mathbb{R}^{N}$ be a closed set with zero
measure, Then there exists $\varphi\in
C^{\infty}_{0}(\mathbb{R}^{N}\setminus K)$ with $\varphi\geqslant0$
such that $\int_{\mathbb{R}^{N}}\nabla u\nabla\varphi dx>0$.
\end{lemma}
{\it Proof.} Since $K\subset\mathbb{R}^{N}$ is closed and $u\neq0$,
we can choose a ball $B\subset\subset\mathbb{R}^{N}\setminus K$, and
a nonnegative function $f\in C^{\infty}_{0}(B)\subset
C^{\infty}_{0}(\mathbb{R}^{N}\setminus K)$ such that
$\int_{\mathbb{R}^{N}} u f dx>0$. Otherwise, we should have that
$u(x)=0$ a.e in $x\in\mathbb{R}^{N}\setminus K$, and it follows from
$|K|=0$ that $u(x)=0$ a.e in $x\in\mathbb{R}^{N}$, which contradicts
$u\neq0$ in $D^{1,2}(\mathbb{R}^{N})$. Then the problem
\begin{equation}\nonumber
\left\{\begin{array}{ll}
 -\Delta v =f, \,\,\,x\in B \\
 v = 0,  \,\,\,x\in \partial B
 \end{array}\right.
\end{equation}
has a nontrivial solution $\tilde{\varphi}\geqslant0$ on $B$ and
$\tilde{\varphi}\in C^{\infty}_0(B)$. Setting
\begin{equation}\nonumber
\varphi=\left\{\begin{array}{ll}
 \tilde{\varphi}, \,\,\,x\in B\\
  0,  \,\,\,x\in\mathbb{R}^{N}\setminus B.
 \end{array}\right.
\end{equation}
Hence,
\begin{equation}\nonumber
\int_{\mathbb{R}^{N}}\nabla u\nabla\varphi dx=\int_{\mathbb{R}^{N}}
u f dx>0
\end{equation}
$\Box$\\
 Based on Lemmas \ref{L:3.4} and \ref{Le:gz}, we prove now the following theorem, which is important for proving our main Theorems \ref{th1.1} and \ref{th1.2}.\\

\begin{theorem}\label{th:3.1} Let $\{u_n\}\subset E$ be nonnegative sequence such that $\|u_n\|_E+\int_{\mathbb{R}^3}\phi_{u_n}u_n^2dx\leqslant C$
and
\begin{equation} \label{eq:3.1.0}
\int_{\mathbb{R}^{3}}[\nabla u_n\nabla
\varphi+(\frac{1}{|y|^{\alpha}}+\lambda_n)u_n\varphi]dx+\int_{\mathbb{R}^{3}}\phi_{u_n}(x)u_n\varphi
dx =\int_{\mathbb{R}^{3}}{u_n}^{p}\varphi dx+o(1),
 \end{equation}
 for any $\varphi\in C^\infty_0(\mathbb{R}^3\setminus T)$, where
 $\alpha\geq0$, $p\in(2,5)$ and $\lambda_n\geqslant0$ with
 $\lambda_n\overset{n}{\rightarrow}\lambda_0<+\infty$. If $\{u_n\}$
 does not converge to 0 in $L^6(\mathbb{R}^3)$, then there exist $\{\tilde{z}_n\}=\{(0,z_n)\}\subset\mathbb{R}^2\times\mathbb{R}$ and nonnegative function $w\in
 E\setminus\{0\}$ such that
$$w_n=T_{1,\tilde{z}_n}u_n\overset{n}{\rightharpoonup} w\text{ weakly
in } E,$$
 and
\begin{equation} \label{eq:3.1.0.0}
\int_{\mathbb{R}^{3}}[\nabla w\nabla
\varphi+(\frac{1}{|y|^{\alpha}}+\lambda_0)w\varphi]dx+\int_{\mathbb{R}^{3}}\phi_{w}(x)w\varphi
dx =\int_{\mathbb{R}^{3}}{w}^{p}\varphi dx,
 \end{equation}
 for any $\varphi\in C^\infty_0(\mathbb{R}^3\setminus T)$. Moreover,
 $\|w\|_E+\int_{\mathbb{R}^3}\phi_{w}w^2dx\leqslant C$ and $w\in C^2(\mathbb{R}^3\setminus T)$.
\end{theorem}
{\it Proof.} If $\{u_n\}\subset E$ does not converges to 0 in
$L^6(\mathbb{R}^3)$, by Lemma \ref{L:3.3} with $N=3$, there exist
$\{s_n\}\subset(0,+\infty)$ and $\{x_n\}\subset\mathbb{R}^3$ with
$x_n=(y_n,z_n)\in\mathbb{R}^2\times\mathbb{R}$ such that
\begin{equation}\label{eq:3.1.1}
 T_{s_n,x_n}u_n\overset{n}{\rightharpoonup}u\neq 0\ \ \text{weakly in
 } L^6(\mathbb{R}^3).
\end{equation}
Let
\begin{equation}\label{eq:3.1.2}
 \tilde{z}_n=(0,z_n)\in \mathbb{R}^2\times\mathbb{R}^1,\ \ w_n=
 T_{1,\tilde{z}_n}u_n=T(1,\tilde{z}_n)u_n(x).
\end{equation}
By (\ref{eq:3.1.1}) and Lemma \ref{L:3.4} with $N=3$, we have that
\begin{equation}\label{eq:3.1.2.0}
 v_n=T_{s_n,0}w_n\overset{n}{\rightharpoonup}v\not\equiv0, \text{ weakly
 in }
 D_s^{1,2}(\mathbb{R}^3),
\end{equation}
where $v$ is nonnegative. And we claim that $s_n>l>0$ for all
$n\in\mathbb{N}$. Indeed, since $-\Delta \phi_{u_n}=u_n^2$, we
easily conclude
$$\int_{\mathbb{R}^3}|u_n|^3dx=\int_{\mathbb{R}^3}\nabla\phi_{u_n}\nabla u_ndx\text{ and }\int_{\mathbb{R}^3}\phi_{u_n}u_n^2dx=\int_{\mathbb{R}^3}|\nabla\phi_{u_n}|^2dx.$$
By using H\"older inequality, we deduce that
$$2\int_{\mathbb{R}^3}|u_n|^3dx\leq\int_{\mathbb{R}^3}|\nabla u_n|^2dx+\int_{\mathbb{R}^3}|\nabla\phi_{u_n}|^2dx=\int_{\mathbb{R}^3}|\nabla u_n|^2dx+\int_{\mathbb{R}^3}\phi_{u_n}u_n^2dx\leq C.$$
So, by using Lemma \ref{L:3.4} with $N=3$ and $q=3$, there is $l>0$
such that $s_n>l$ for all
$n\in\mathbb{N}$.\\
 \indent{\bf Step1:} There exists $L>l>0$
such that $s_n<L$ for
 $n\in\mathbb{N}$ large.\\
 Recalling the definition of $T$ in (\ref{T}), we have $|T|=0$. Since
 $v\geqslant0$, by Lemma \ref{Le:gz}, we have a nonnegative function $\varphi_1\in C^\infty_0(\mathbb{R}^3\setminus
 T)$ such that
\begin{equation}\nonumber
\int_{\mathbb{R}^3}\nabla v\nabla\varphi_1 dx>0.
\end{equation}
It follows from (\ref{eq:3.1.2}) and (\ref{eq:3.1.2.0}) that
\begin{equation}\label{eq:3.1.3}
\int_{\mathbb{R}^{3}}(\nabla(
 T_{s_n,\tilde{z}_n}u_{n})\nabla \varphi_1
dx\rightarrow\int_{\mathbb{R}^{3}}\nabla v\nabla \varphi_1 dx>0.
\end{equation}
 Noting that
$T^{-1}_{s_n,\tilde{z}_n}\varphi_1(x)=s_{n}^{\frac{1}{2}}\varphi_1(s_n
x-s_n\tilde{z}_n)$, then $T^{-1}_{s_n,\tilde{z}_n}\varphi_1(x)\in
C^\infty_0(\mathbb{R}^3\setminus T)$, by (\ref{eq:3.1.0}), as
$n\rightarrow+\infty$, we have that

\begin{eqnarray}\nonumber
\int_{\mathbb{R}^{3}}\phi_{u_{n}}u_{n}T^{-1}_{s_n,\tilde{z}_n}\varphi_1
dx&+&\int_{\mathbb{R}^{3}}[\nabla u_{n}\nabla
(T^{-1}_{s_n,\tilde{z}_n}\varphi_1)+(\lambda_n+\frac{1}{|y|^\alpha})u_{n}T^{-1}_{s_n,\tilde{z}_n}\varphi_1]dx
\nonumber\\
&=&\int_{\mathbb{R}^{3}}u_{n}^{p}T^{-1}_{s_n,\tilde{z}_n}\varphi_1
dx+o(1)\nonumber.
\end{eqnarray}
It follows from $u_n\geqslant0$ and $\lambda_n\geqslant0$ that
\begin{eqnarray}\nonumber
\int_{\mathbb{R}^{3}}\nabla u_{n}\nabla
(T^{-1}_{s_n,\tilde{z}_n}\varphi_1)dx\leqslant\int_{\mathbb{R}^{3}}u_{n}^{p}T^{-1}_{s_n,\tilde{z}_n}\varphi_1
dx+o(1)\nonumber.
\end{eqnarray}
That is
\begin{eqnarray}\nonumber
\int_{\mathbb{R}^{3}}(\nabla( T_{s_n,\tilde{z}_n}u_{n})\nabla
\varphi_1 dx &\leqslant&s_n^{\frac{p-5}{2}}\int_{\mathbb{R}^{3}}{(
T_{s_n,\tilde{z}_n}u_{n})}^{p}\varphi_1
dx+o(1)\nonumber\\
&\leqslant&Cs_n^{\frac{p-5}{2}}\int_{supp
\varphi_1}{(T_{s_n,\tilde{z}_n}u_{n})}^{p} dx+o(1)\nonumber\\
&\leqslant&Cs_n^{\frac{p-5}{2}}{\|
T_{s_n,\tilde{z}_n}u_{n}\|}_{6}^{p}
+o(1)\,\,\text{for}\,\,2<p<5\nonumber\\
&\leqslant&Cs_n^{\frac{p-5}{2}}{\|\nabla u_{n}\|}_{2}^{p}
+o(1)\quad\quad\text{by (\ref{eq:3.1})}\nonumber.
\end{eqnarray}
Since $\{u_n\}$ is bounded in $E$, if $s_n\rightarrow\infty$, it
follows that
$\limsup\limits_{n\rightarrow\infty}\int_{\mathbb{R}^{3}}(\nabla(
 T_{s_n,\tilde{z}_n}u_{n})\nabla \varphi
dx\leqslant0$, which contradicts with (\ref{eq:3.1.3}).\\
\indent{\bf Step 2:} $\{w_n\}$ is a bounded sequence in $E$ such
that for any $\varphi\in C^\infty_0(\mathbb{R}^3\setminus T)$, as
$n\rightarrow+\infty$,
\begin{equation} \label{eq:3.1.4}
\int_{\mathbb{R}^{3}}[\nabla w_n\nabla
\varphi+(\frac{1}{|y|^{\alpha}}+\lambda_n)w_n\varphi]dx+\int_{\mathbb{R}^{3}}\phi_{w_n}(x)w_n\varphi
dx =\int_{\mathbb{R}^{3}}{w_n}^{p}\varphi dx+o(1).
 \end{equation}
By the definition of $T_{s,\xi}$ in (\ref{eq:3.2}), we have
\begin{equation}\nonumber
 \|\nabla(T_{1,\tilde{z}_n}u_n)\|_2=\|\nabla u_n\|_2,\ \
 \  \int_{\mathbb{R}^{3}}\frac{|T_{1,\tilde{z}_n}u_n|^2}{|y|^\alpha}dx
 =\int_{\mathbb{R}^{3}}\frac{|u_n|^2}{|y|^\alpha}dx,
\end{equation}
hence,
$\|w_n\|_{E}^{2}=\|T_{1,\tilde{z}_n}u_n\|_{E}^{2}=\|u_n\|_{E}^{2}$
and $\{w_n\}$ is bounded in $E$. By the definitions of
$T_{1,\tilde{z}_n}$ in (\ref{eq:3.2}) and $\phi_u$ in
(\ref{eq:1.7}), it is easy to see that
\begin{eqnarray} \nonumber
\int_{\mathbb{R}^{3}}\phi_{u_n}u_nT^{-1}_{1,\tilde{z}_n}\varphi dx
=\int_{\mathbb{R}^{3}}T_{1,\tilde{z}_n}(\phi_{u_n}u_n)\varphi dx
=\int_{\mathbb{R}^{3}}\phi_{w_n}w_n\varphi dx,
 \end{eqnarray}
\begin{equation} \nonumber
\int_{\mathbb{R}^{3}}[\nabla u_n\nabla
T^{-1}_{1,\tilde{z}_n}\varphi+(\frac{1}{|y|^{\alpha}}+\lambda_n)u_nT^{-1}_{1,\tilde{z}_n}\varphi]dx=\int_{\mathbb{R}^{3}}[\nabla
w_n\nabla \varphi+(\frac{1}{|y|^{\alpha}}+\lambda_n)w_n\varphi]dx,
 \end{equation}
 and
\begin{equation} \nonumber
\int_{\mathbb{R}^{3}}u_n^{p}T^{-1}_{1,\tilde{z}_n}\varphi
dx=\int_{\mathbb{R}^{3}}w_n^{p}\varphi dx.
 \end{equation}
 Those and (\ref{eq:3.1.0}) imply that (\ref{eq:3.1.4}) holds.\\
\indent{\bf Step 3:} $w_n\overset{n}{\rightharpoonup}w\not\equiv0$
in $E$ and $w(x)\geqslant0$ a.e. in $x\in \mathbb{R}^3$.\\
By Step 1, there exists $s_0\in[l,L]$ such that, passing to
subsequence, $s_n\overset{n}{\rightarrow}s_0$. Then, it follows from
(\ref{eq:3.1.2.0}) and Lemma \ref{L:3.2.0} that
\begin{equation}\label{eq:3.1.5}
w_n=T^{-1}_{s_n,0}v_n\overset{n}{\rightharpoonup}T_{\frac{1}{s_0},0}v\not\equiv0\
 \ \ \text{weakly in } D_s^{1,2}(\mathbb{R}^{3}).
\end{equation}
By Step 2, there exists $w\in E$ such that, passing to a
subsequence, $w_n\overset{n}{\rightharpoonup}w$ weakly in $E$, since
$E\subset D_s^{1,2}(\mathbb{R}^3)$, we have
$(D_s^{1,2}(\mathbb{R}^3))^\ast\subset E^\ast$, hence
$w_n\overset{n}{\rightharpoonup}w$ weakly in $
D_s^{1,2}(\mathbb{R}^3)$, it follows from (\ref{eq:3.1.5}) that
$w=T_{\frac{1}{s_0},0}v\not\equiv0$ and $w(x)\geqslant0$ a.e. in
$x\in\mathbb{R}^3$, since $v\geq0$ in (\ref{eq:3.1.2.0}).\\
\indent{\bf Step 4:} $\phi_w\in D^{1,2}(\mathbb{R}^3)$ and
(\ref{eq:3.1.0.0}) holds. \\
For each $n\in\mathbb{N}$, $\|\nabla
\phi_{w_n}\|_2^2=\int_{\mathbb{R}^3}\phi_{w_n}w_n^2dx=\int_{\mathbb{R}^3}\phi_{u_n}u_n^2dx$,
hence, $\int_{\mathbb{R}^3}\phi_{u_n}u_n^2dx<C$ implies that
$\{\phi_{w_n}\}$ is bounded in $D_s^{1,2}(\mathbb{R}^3)$. So, there
exists $\phi\in D_s^{1,2}(\mathbb{R}^3)$ such that
$\phi_{w_n}\overset{n}{\rightharpoonup}\phi$ weakly in
$D_s^{1,2}(\mathbb{R}^3)$, that is
\begin{equation}\label{eq:3.1.6}
\int_{\mathbb{R}^3}\nabla \phi_{w_n}\nabla\varphi
dx\overset{n}{\rightarrow} \int_{\mathbb{R}^3}\nabla
\phi\nabla\varphi dx, \text{ for any }\varphi\in
C^\infty_0(\mathbb{R}^3).
\end{equation}
On the other hand, for any $\varphi\in C^\infty_0(\mathbb{R}^3)$, we
have
\begin{equation}\label{eq:3.1.7}
\int_{\mathbb{R}^3}\nabla \phi_{w_n}\nabla\varphi
dx=\int_{\mathbb{R}^3} w_n^2\varphi dx\ \text{ and }
\int_{\mathbb{R}^3}w_n^2\varphi dx \overset{n}{\rightarrow}
\int_{\mathbb{R}^3}w^2\varphi dx.
\end{equation}
It follows from (\ref{eq:3.1.6}) and (\ref{eq:3.1.7}) that
\begin{equation}\nonumber
\int_{\mathbb{R}^3}\nabla \phi\nabla\varphi dx=
\int_{\mathbb{R}^3}w^2\varphi dx\text{ for any }\varphi\in
C^\infty_0(\mathbb{R}^3).
\end{equation}
So, $\phi$ is a solution of $-\Delta\phi=w^2$ in the sense of
distribution. Since $w\in E\subset L^6(\mathbb{R}^3)$,
$\phi_w(x)=\int_{\mathbb{R}^3}\frac{w^2(y)}{|x-y|}dy\in
W^{2,3}(\mathbb{R}^3)$ by Theorem 9.9 in
\cite{DavidGilbarg-Neil.S.Trudinger}, hence $\phi_w$ satisfies
$-\Delta\phi_w=w^2$ in the sense of distribution (Theorem 6.21 in
\cite{EHLiebMLoss-analysis}). By uniqueness, we have $\phi_w=\phi\in
D_s^{1,2}(\mathbb{R}^3)$. It follows from (\ref{eq:3.1.6}) that
\begin{equation}\nonumber
\phi_{w_n}\overset{n}{\rightharpoonup}\phi_w\  \text{ weakly in }\ \
\ D_s^{1,2}(\mathbb{R}^3).
\end{equation}
Then (see (3.18) in \cite{JZ-ActMS} for the details), for any
$\varphi\in C^\infty_0(\mathbb{R}^3)$, we have
\begin{equation}\nonumber
\int_{\mathbb{R}^3}\phi_{w_n}(x)w_n\varphi
dx\overset{n}{\rightarrow}\int_{\mathbb{R}^3}\phi_{w}(x)w\varphi dx.
\end{equation}
For each bounded domain $\Omega\subset\mathbb{R}^3$ and $q\in(1,6)$,
it follows from (\ref{eq:3.1.6}) and the compactness of Sobolev
embedding that $w_n\overset{n}{\rightarrow}w$ strongly in
$L^q(\Omega)$. Hence, for any $\varphi\in
C^\infty_0(\mathbb{R}^3\setminus T)$,
\begin{equation}\nonumber
\int_{\mathbb{R}^{3}}(\nabla w_n\nabla
\varphi+(\frac{1}{|y|^{\alpha}}+\lambda_n)w_n
\varphi)dx\overset{n}{\longrightarrow}\int_{\mathbb{R}^{3}}(\nabla
w\nabla\varphi+(\frac{1}{|y|^{\alpha}}+\lambda_0)w\varphi)dx
\end{equation}
and
\begin{equation}\nonumber
\int_{\mathbb{R}^{3}}w_n^{p} \varphi
dx\overset{n}{\longrightarrow}\int_{\mathbb{R}^{3}}w^{p}\varphi dx.
\end{equation}
Those and (\ref{eq:3.1.4}) imply that (\ref{eq:3.1.0.0}) holds.\\
\indent{\bf Step 5.} $\| w\|_{E}+\int_{\mathbb{R}^{3}}\phi_{w}w^2
dx<C$.\\
By Step 3, we have $w_n\overset{n}{\rightharpoonup}w \text{ weakly
in } E$, and Step 4 implies that
\begin{equation}\nonumber
\int_{\mathbb{R}^{3}}\phi_{w}w^2 dx= \|\nabla \phi_w\|_2^2\ \
\text{and } \phi_{w_n}\overset{n}{\rightharpoonup}\phi_w\ \text{
weakly in }\ D^{1,2}(\mathbb{R}^{3}),
\end{equation}
and the lower semi-continuity of norm implies that
\begin{equation}\nonumber
\|w\|_E\leqslant\liminf\limits_{n\rightarrow+\infty}\|w_n\|_E,
\end{equation}
and
\begin{equation}\nonumber
\int_{\mathbb{R}^{3}}\phi_{w}w^2 dx=\|\nabla
\phi_{w}\|_2^2\leqslant\liminf\limits_{n\rightarrow+\infty}\|\nabla
\phi_{w_n}\|_2^2=\liminf\limits_{n\rightarrow+\infty}\int_{\mathbb{R}^{3}}\phi_{w_n}w_n^2
dx.
\end{equation}
Hence, by (\ref{eq:3.1.2}), we have
\begin{eqnarray}
\|w\|_E+\int_{\mathbb{R}^{3}}\phi_{w}w^2
dx&\leqslant&\liminf\limits_{n\rightarrow+\infty}\left\{\|w_n\|_E+
\int_{\mathbb{R}^{3}}\phi_{w_n}w_n^2
dx\right\}\nonumber\\
&=&\liminf\limits_{n\rightarrow+\infty}\left\{\|u_n\|_E+
\int_{\mathbb{R}^{3}}\phi_{u_n}u_n^2 dx\right\}\leqslant C.
\nonumber
\end{eqnarray}
\indent{\bf Step 6.} $w(x)\in C^2(\mathbb{R}^3\setminus T)$.\\
Since $\lambda_0\geq0$ and $w(x)\geq0$ for a.e. $x\in\mathbb{R}^3$, it follows from (\ref{eq:3.1.0.0}) that,  for any nonnegative function $v\in
C^\infty(\mathbb{R}^3\setminus T)$,
\begin{equation} \label{inf}
\int_{\mathbb{R}^{3}}\nabla w\nabla vdx
\leq\int_{\mathbb{R}^{3}}{w}^{p}v dx.
 \end{equation}
Then, Lemma \ref{L:4.1} in section 4 implies that (\ref{inf}) holds also
for any nonnegative function $v\in H^1(\mathbb{R}^3)$. Note that,
for any nonnegative function $\varphi\in C^\infty_0(\mathbb{R}^3)$
and any nonnegative piecewise smooth function $h$ on $[0,+\infty)$, $h(w)\varphi\in H^1(\mathbb{R}^3)$. Take
$v=h(w)\varphi$ in (\ref{inf}), then we see that
(\ref{eq:4.3.1}) in section 4 holds with $u=w$ and $N=3$. Hence, by Lemma
\ref{Le:4.3}, we have $w\in L^\infty(\mathbb{R}^3)$.
 Let $\Omega\subset\subset\mathbb{R}^3\setminus T$ be a bounded
 domain with smooth boundary, then $\frac{1}{|y|}$ is a smooth
function in $\Omega$ and $w\in W^{1,2}(\Omega)$ is a weak solution
of
\begin{equation}\label{eq:regu}
-\Delta w(x)=f(x),\ \ \ x\in\Omega,
\end{equation}
where
$f(x)=|w|^{p-1}w(x)-\phi_w(x)w(x)-(\lambda_0+\frac{1}{|y|})w(x)$.
Since $w, \phi_w\in W^{1,2}(\Omega)$ and $w\in L^\infty(\Omega)$, we have $f(x)\in
W^{1,2}(\Omega)$. By using Theorem 8.10 in
\cite{DavidGilbarg-Neil.S.Trudinger}, we get $w\in
W_{loc}^{3,2}(\Omega)$. Then, Sobolev imbedding theorem implies that
$w\in C_{loc}^{1/4}(\Omega)$, hence $\phi_w(x)\in
C_{loc}^{2,1/4}(\Omega)$ since $\phi_w(x)$ is a weak solution of
$-\Delta \phi(x)=w^2(x)$ in $D^{1,2}(\Omega)$. It follows that
$f(x)\in C_{loc}^{1/4}(\Omega)$. By applying Theorem 9.19 in
\cite{DavidGilbarg-Neil.S.Trudinger} to (\ref{eq:regu}), we have
$w\in C_{loc}^{2,1/4}(\Omega)$. So $w\in C^2(\mathbb{R}^3\setminus
T)$.\ \ \ \ \ \ \ \ \ \ \ \ \ \  $\Box$

{\bf Proof of Theorem \ref{th1.1}}
 Let $\{u_n\}\subset H$ be the
bounded nonnegative (PS) sequence obtained by Lemma \ref{L:2.6},
then there exists $C>0$, which is independent of $\lambda$ if
$\lambda\in(0,1)$, such that
\begin{equation}\label{eq:th1.1.0}
\|u_n\|^2_{H}+\int_{\mathbb{R}^{3}}\phi_{u_n}{u_n}^2 dx\leqslant
C\text{ and } \ u_n(x)\geqslant0\ \text{a.e. in}\ \ x\in
\mathbb{R}^3.
\end{equation}
Hence,
\begin{equation}\nonumber
\|u_n\|^2_{E}+\int_{\mathbb{R}^{3}}\phi_{u_n}{u_n}^2 dx\leqslant C,\
\ u_n(x)\geqslant0\ \text{a.e. in}\ \ x\in \mathbb{R}^3.
\end{equation}
 And (\ref{eq:2.10}) implies that (\ref{eq:3.1.0})
holds with $\lambda_n\equiv\lambda>0$. If $\{u_n\}$ does not
converges to 0 in $L^6(\mathbb{R}^3)$, by Theorem \ref{th:3.1},
there exist
$\{\tilde{z}_n\}=\{(0,z_n)\}\subset\mathbb{R}^2\times\mathbb{R}$ and
nonnegative function $w\in
 E\setminus\{0\}$ such that
 \begin{equation} \label{eq:th1.1.1}
w_n=T_{1,\tilde{z}_n}u_n\overset{n}{\rightharpoonup} w\text{ weakly
in } E,
 \end{equation}
\begin{equation}\nonumber
\int_{\mathbb{R}^{3}}[\nabla w\nabla
\varphi+(\frac{1}{|y|^{\alpha}}+\lambda)w\varphi]dx+\int_{\mathbb{R}^{3}}\phi_{w}(x)w\varphi
dx =\int_{\mathbb{R}^{3}}{w}^{p}\varphi dx,
 \end{equation}
 for any $\varphi\in C^\infty_0(\mathbb{R}^3\setminus T)$, i.e., $w$ is a weak solution of (\ref{eq:1.1}) in $E$. Moreover, $w\in C^2(\mathbb{R}^3\setminus T)$ and
\begin{equation}\label{eq:th1.1.3}
\|w\|_E+\int_{\mathbb{R}^3}\phi_{w}w^2dx\leqslant
 C.
 \end{equation}

Now, we claim that $w\in H$. In fact,
by (\ref{eq:3.2}) and (\ref{eq:th1.1.1}), we have
$\|w_n\|_H=\|u_n\|_H$ and $\|w_n\|_H$ is bounded, so there exists
$w^\ast\in H$ such that
\begin{equation}\label{eq:th1.1.2}
w_n\overset{n}{\rightharpoonup}w^\ast \text{ weakly in } H\text{ and
} w_n(x)\overset{n}{\rightarrow}w^\ast(x),\ \text{a.e. in
}x\in\mathbb{R}^3.
\end{equation}
On the other hand, (\ref{eq:th1.1.1}) implies that
\begin{equation}\nonumber
 w_n(x)\overset{n}{\rightarrow}w(x),\ \text{a.e. in
}x\in\mathbb{R}^3.
\end{equation}
This and (\ref{eq:th1.1.2}) show that $w=w^\ast\in H$. Moreover, if
$\lambda\in(0,1)$, Lemma \ref{L:2.6} shows that there exists $M>0$
independent of $\lambda\in(0,1)$ such that (\ref{eq:th1.1.0}) holds
with $C=M$, then (\ref{eq:th1.1.3}) holds with $C=M$. Hence, to
complete the proof of Theorem \ref{th1.1}, we only need to prove
that $\{u_n\}$ cannot converges to 0 in $L^6(\mathbb{R}^3)$. For
$r\in(2,6)$, by H\"{o}lder inequality we have
\begin{equation*}
 \int_{\mathbb{R}^3}{|u_n|}^rdx=\int_{\mathbb{R}^3}{|u_n|}^{\frac{2}{q}}{|u_n|}^{\frac{6}{q{'}}}dx
 \leqslant\|u_n\|_{2}^{\frac{2}{q}}\|u_n\|_{6}^{\frac{6}{q{'}}}
\end{equation*}
where $q=\frac{4}{6-r}>1$, $q{'}=\frac{q}{q-1}=\frac{4}{r-2}>1$.
Hence, if $u_n\overset{n}{\rightarrow} 0$ in $L^{6}(\mathbb{R}^3)$,
then $u_n\overset{n}{\rightarrow} 0$ in $L^{r}(\mathbb{R}^3)$ for
$r\in(2,6)$, this and (\ref{eq:1.8}) imply that
$\int_{\mathbb{R}^3}\phi_{u_n}(x){u_n}^2dx\overset{n}{\rightarrow}
0$. Therefore, by (\ref{eq:2.10}) we have that, for $p\in(2,5)$,
\begin{eqnarray*}\nonumber
c&=&\lim\limits_{n\rightarrow\infty}\left[I(u_n)-\frac{1}{2}I{'}(u_n)u_n\right]\\
&=&\lim\limits_{n\rightarrow\infty}\left[-\frac{1}{4}\int_{\mathbb{R}^3}\phi_{u_n}(x){u_n}^2dx+\frac{p-3}{2(p+1)}\int_{\mathbb{R}^3}|u_n|^{p+1}dx\right]=0,
\end{eqnarray*}
this is impossible since $c>0$.  $\Box$\\
\section{Existence for $\lambda=0$: Proof of Theorem \ref{th1.2}.}
We need more lemmas as follows to prove Theorem \ref{th1.2}. For any
$N\geqslant3$ and domain $\Omega\subset \mathbb{R}^N$($\Omega$ can
be bounded or unbounded), let $\Gamma\subset\Omega$ be a closed
Manifold with $codim \Gamma=k\geqslant2$. Then,
\begin{lemma}\label{L:4.0}
$C^{\infty}_0(\Omega\setminus\Gamma)$ is dense in $H^1_0(\Omega)$.
\end{lemma}
{\bf Proof:}
 For each $u\in H_0^1(\Omega)\cap
{C_0^\infty(\Omega\setminus\Gamma)}^{\bot}$ and $\tilde{\varphi}\in
C_0^\infty(\Omega\setminus\Gamma)$, we have
\begin{eqnarray}\label{1}
\langle u,\tilde{\varphi}\rangle_{H_0^1(\Omega)}=0,
\end{eqnarray}
since $C_0^\infty(\Omega\setminus\Gamma)$ is dense in
$H_0^1(\Omega\setminus\Gamma)$, it follows that
\begin{eqnarray}\label{2}
\langle u,\psi\rangle_{H_0^1(\Omega)}=0 \text{ for any } \psi\in
H_0^1(\Omega\setminus\Gamma).
\end{eqnarray}
It is true that $ C_0^\infty(\Omega\setminus\Gamma)$ is dense in
$H_0^1(\Omega)$ if $C_0^\infty(\Omega\setminus\Gamma)^{\bot}\cap
H_0^1(\Omega)=\{0\}$. Hence, we only need to show that (\ref{1})
holds with all
$\tilde{\varphi}\in C_0^\infty(\Omega)$ as follows.\\
For any $\varphi\in C^\infty_0(\Omega)$,
 let $\Omega_0=\text{supp}\varphi$.
If $\Omega_0\cap\Gamma=\emptyset$, then $\varphi\in
C^\infty_0(\Omega\setminus\Gamma)$ and (\ref{1}) holds with
$\tilde{\varphi}=\varphi$. Otherwise
$\Omega_0\cap\Gamma\neq\emptyset$, setting
$\Gamma_0=\Omega_0\cap\Gamma$, and for any $d>0$ small enough that
we have set $\Gamma_d:=\{x\in \Omega:dist(x,\Gamma_0)<d\}\subset
\Omega$. Let
\begin{equation}\nonumber
\psi_d(x):=\left\{\begin{array}{ll}
 \frac{dist(x,\Gamma_{2d})}{d}, \  \hfill x\in \Gamma_{3d}, \\
 1,\   \hfill  x\in\Omega\setminus\Gamma_{3d},
 \end{array}\right.
\end{equation}
then $\psi_d(x)\in C^{0,1}(\Omega) \text{ and }\|\psi_d\|_{
C^{0,1}(\Omega)}\leqslant\frac{1}{d}$. Let
$\varphi_d:=\varphi(1-\psi_d)$, we have $ \varphi\psi_d\in
H_0^1(\Omega\setminus \Gamma )\text{ and }  \ \varphi_d\in
H_0^1(\Gamma_{3d})$. It follows from (\ref{2}) that
\begin{eqnarray}\label{4}
\langle u,\varphi\rangle_{H^1}&=&\langle
u,\varphi_d+\varphi\psi_d\rangle_{H^1}
=\langle u,\varphi_d\rangle_{H^1}+\langle u,\varphi\psi_d\rangle_{H^1}\nonumber\\
&=&\langle u,\varphi_d\rangle_{H^1}
\leqslant\|u\|_{H^1(\Gamma_{3d})}\|\varphi_d\|_{H^1(\Gamma_{3d})}.
\end{eqnarray}
By the definition of $\varphi_d$,  we have
\begin{eqnarray}\label{5}
\|\varphi_d\|^2_{ L^2(\Gamma_{3d})}=\int_{\Gamma_{3d}}\varphi^2_d
dx\leqslant
4\|\varphi\|^2_{L^{\infty}(\Omega)}|\Gamma_{3d}|\overset{d\rightarrow
0}{\longrightarrow}0,
\end{eqnarray}
\begin{eqnarray}\label{6}
\|\nabla \varphi_d\|^2_{
L^2(\Gamma_{3d})}&=&\int_{\Gamma_{3d}}|\nabla \varphi_d|^2
dx\leqslant C
\|\varphi\|^2_{C^{1}(\Omega)}|\Gamma_{3d}|(1+\frac{1}{d^2})\nonumber\\
\text{ since } k\geqslant2&&\overset{codim \Gamma=k}{\leqslant}C
d^{k-2}\leqslant C.
\end{eqnarray}
And $|\Gamma_{3d}|\overset{d\rightarrow0}{\longrightarrow0}$ implies
that
\begin{eqnarray}\label{7}
\|u\|_{H^1(\Gamma_{3d})}\overset{d\rightarrow 0}{\longrightarrow}0.
\end{eqnarray}
It follows from (\ref{4}) to (\ref{7}) that (\ref{1}) holds for all $\tilde{\varphi}\in C^\infty_0(\Omega)$. $\Box$\\
\begin{lemma}\label{L:4.1}
$\{\varphi\in H^{1}_0(\Omega\setminus\Gamma):\varphi(x)\geqslant0\}$
is dense in $\{\varphi\in H^1_0(\Omega):\varphi(x)\geqslant0\}$.
\end{lemma}
 {\bf Proof:} Lemma \ref{L:4.0} shows that for any $u(x)\in
H_0^{1}(\Omega)$, there exist $\{\varphi_n(x)\}\subset
C^{\infty}_0(\Omega\setminus\Gamma)$ such that
\begin{eqnarray}\nonumber
\|\varphi_n-u\|_{H^1(\Omega)}\overset{n}{\rightarrow}0.
\end{eqnarray}
 This lemma is proved if we have
\begin{eqnarray}\nonumber
\left\||\varphi_n|-|u|\right\|_{H^1(\Omega)}\overset{n}{\rightarrow}0,
\end{eqnarray}
which is true by the following two facts,
\begin{eqnarray*}
0&\leqslant&\left\||\varphi_n|-|u|\right\|^2_{2}=\int_{\Omega}|\varphi_n|^2+|u|^2-2|\varphi_n||u|dx\\
&\leqslant& \int_{\Omega}\varphi_n^2+u^2-2\varphi_nudx
=\|\varphi_n-u\|^2_{2}\overset{n}{\rightarrow}0.
\end{eqnarray*}
\begin{eqnarray*}
0&\leqslant&\left\||\varphi_n|-|u|\right\|^2_{D^{1,2}}=\int_{\Omega}|\nabla|\varphi_n||^2+|\nabla|u||^2-2\nabla|\varphi_n|\nabla|u|dx\\
&=&\int_{\Omega}|\nabla(\varphi_n-
u)|^2dx+4\int_{\Omega}\nabla\varphi^{+}_n\nabla u^{-}
+\nabla\varphi^{-}_n\nabla u^{+}dx
 \overset{n}{\longrightarrow}0.  \ \ \Box
\end{eqnarray*}
\begin{lemma}\label{Le:4.3}(Lemma 3.2 of \cite{JZ-reprint})
Let $N\geqslant3$, $p\in(1,\frac{N+2}{N-2})$ and let $u\in
D^{1,2}(\mathbb{R}^N)\setminus\{0\}$ be a nonnegative function such
that
\begin{equation}\label{eq:4.3.1}
\int_{\mathbb{R}^{N}}\nabla u\nabla (h(u)\varphi)
dx\leqslant\int_{\mathbb{R}^{N}}|u|^{p-1}uh(u)\varphi dx,
\end{equation}
holds for any nonnegative $\varphi\in C_0^\infty(\mathbb{R}^N)$ and
any nonnegative piecewise smooth function $h$ on $[0,+\infty)$ with
$h'\in L^\infty(\mathbb{R})$. Then, $u\in L^\infty(\mathbb{R}^N)$
and there exist $C_1>0$ and $C_2>0$, which depend only on $N$ and
$p$, such that
\begin{eqnarray}\nonumber
\|u\|_{\infty}\leqslant
C_1\left(1+\|u\|_{2^\ast}^{C_2}\right)\|u\|_{2^\ast}\nonumber.
\end{eqnarray}
\end{lemma}

\begin{lemma}\label{Le:4.4.0}
For $p>2$, let $(u,\phi)\in H^1(\mathbb{R}^3)\cap
L^{p+1}(\mathbb{R}^3)\times {D}^{1,2}(\mathbb{R}^3)$ be a nontrivial
nonnegative weak solution of the following problem
\begin{equation}\label{eq:est1}
\left\{\begin{array}{ll}
 -\Delta u +\mu\phi (x) u\leq|u|^{p-1}u,\,\,\, x\in\mathbb{R}^3, \\
 -\Delta\phi = u^2,\,\, x\in \mathbb{R}^3,\\
 \end{array}\right.
\end{equation}
 where $\mu>0$.
 Then
\begin{equation}\nonumber
\|u\|_{\infty}>\mu^{\frac{1}{2(p-2)}}.
\end{equation}
\end{lemma}
{\it Proof:} By assumption, $(u,\phi)\in
 H^1(\mathbb{R}^3)\cap L^{p+1}(\mathbb{R}^3)\times {D}^{1,2}(\mathbb{R}^3)$  is a
 weak solution of
(\ref{eq:est1}), then, for any nonnegative function $v\in
H^1(\mathbb{R}^3)\cap L^{p+1}(\mathbb{R}^3)$, we have
\begin{equation}\label{eq:est2}
\int_{\mathbb{R}^3} \nabla u \nabla
vdx+\mu\int_{\mathbb{R}^3}\phi(x)uvdx
 -\int_{\mathbb{R}^3}|u|^{p-1}u vdx\leq0,
\end{equation}
\begin{equation}\label{eq:est3}
\int_{\mathbb{R}^3} \nabla \phi \nabla vdx
 =\int_{\mathbb{R}^3} u^{2} vdx.
\end{equation}
For $c>0$, adding $c\int_{\mathbb{R}^3}u^2vdx$ to both sides of
(\ref{eq:est2}), and using (\ref{eq:est3}) we get that
\begin{equation}\label{eq:est3.0}
\begin{split}
\int_{\mathbb{R}^3} \nabla u \nabla
vdx&+\int_{\mathbb{R}^3}[cu^2-|u|^{p-1}u]v
dx+\mu\int_{\mathbb{R}^3}\phi (x)uvdx\\
&\leq c\int_{\mathbb{R}^3} \nabla \phi\nabla vdx,\ \text{ for any }
v\in
 H^1(\mathbb{R}^3)\cap L^{p+1}(\mathbb{R}^3).
 \end{split}
\end{equation}
In the following, we mean that $w_+(x)=\max\{0,w(x)\}$ for any function
$w(x)$ on $\mathbb{R}^3$. For the above $c>0$, we let $\epsilon>0$
small,
\begin{equation}\label{eq:est3.1}
w_1(x)=(u(x)-c\phi(x)-\epsilon)^+ \text{ and }
\Omega_1=\{x\in\Omega:w_1(x)>0\}.
\end{equation}
It is easy to see that
$u(x)\overset{|x|\rightarrow+\infty}{\rightarrow}0$ and
$\phi(x)\geq0$ a.e. $x\in \mathbb{R}^3$, then $w_1\in
H^1(\mathbb{R}^3)\cap L^{p+1}(\mathbb{R}^3)$ and $u(x)|_{\Omega_1}>
c\phi(x)>0$. Taking $v(x)=w_1(x)$ in (\ref{eq:est3.0}), we see that
\begin{equation}\label{eq:est4}
\int_{\Omega_1} \nabla u \nabla w_1 dx+\int_{\Omega_1}[c
u^2-|u|^{p-1}u]w_1dx\leqslant
 c\int_{\Omega_1} \nabla \phi\nabla w_1 dx.
\end{equation}
However, for all $x\in\Omega_1$ we have $cu^2-|u|^{p-1}u\geq0$ if
$c=\delta^{p-2}$ with $\delta=\|u\|_{\infty}$. Then, let
$c=\delta^{p-2}$
 and (\ref{eq:est4}) implies that
$$\int_{\Omega_1} \nabla u \nabla w_1 dx-
 c\int_{\Omega_1} \nabla \phi\nabla w_1 dx\leqslant0,$$
 that is,
\begin{equation}\label{eq:est5}
 \int_{\Omega_1}\nabla (u-\delta^{p-2}\phi)\nabla w_1dx=
 \int_{\Omega_1} |\nabla w_1|^2dx=
 0.
\end{equation}
Hence, either $|\Omega_1|=0$ or
$w_1|_{\Omega_1}\equiv\text{constant}$,
 this means that
$u(x)\leqslant \delta^{p-2}\phi(x)+\epsilon$ a.e.
$x\in\mathbb{R}^3$. Let $\epsilon\rightarrow0$ we have
\begin{equation}\label{eq:est6}
u(x)\leqslant \delta^{p-2}\phi(x),\text{ a.e. in }x\in\mathbb{R}^3,
\end{equation}
 To prove that $\|u\|_{\infty}>\mu^{\frac{1}{2(p-2)}}$, we let
 $v=u$ in (\ref{eq:est2}), it follows that
\begin{equation}\nonumber
\int_{\mathbb{R}^3} |\nabla u
|^2dx+\mu\int_{\mathbb{R}^3}\phi(x)u^2dx
 -\int_{\mathbb{R}^3}u^{p+1}dx\leq0,
\end{equation}
that is,
\begin{equation}\nonumber
\mu\int_{\mathbb{R}^3}\phi(x)|u|^2dx
\leq\int_{\mathbb{R}^3}|u|^{p+1}dx.
\end{equation}
This and (\ref{eq:est6}) show that
\begin{equation}\nonumber
\int_{\mathbb{R}^3}(u^{p-2}-\mu\delta^{2-p})u^3dx\geq0.
\end{equation}
Hence, $\delta^{p-2}\geq\mu\delta^{2-p}$ by $p>2$. On the other
hand, by using $u\neq0$ we have $\delta>0$. Then
$\|u\|_\infty=\delta\geq\mu^{\frac{1}{2(p-2)}}$.
  $\Box$\\

{\bf Proof of Theorem \ref{th1.2}.} By Theorem
\ref{th1.1}, we know that, for each $\lambda\in(0,1)$, problem (\ref{eq:1.1}) has nonnegative
solution $u_\lambda\in H\setminus\{0\}$ such that
$\|u_\lambda\|_E+\int_{\mathbb{R}^3}\phi_{u_\lambda}u_\lambda^2dx\leqslant
M$ and (\ref{eq:3.1.0}) holds with $u_n=u_\lambda$ and
$\lambda_n=\lambda$. Since $u_\lambda\geqslant0$, it follows from (\ref{eq:3.1.0}) that
\begin{equation}\nonumber
\int_{\mathbb{R}^{3}}\nabla u_\lambda\nabla \varphi
dx+\int_{\mathbb{R}^3}\phi_{u_\lambda}(x)u_\lambda\varphi
dx\leqslant\int_{\mathbb{R}^{3}}u_\lambda^p\varphi dx\ \text{ for
all } \varphi\in C^\infty_0(\mathbb{R}^3\setminus T),
\varphi\geqslant0.
\end{equation}
This and Lemma \ref{L:4.1} show that
\begin{equation}\label{eq:th1.2.1}
\int_{\mathbb{R}^{3}}\nabla u_\lambda\nabla v
dx+\int_{\mathbb{R}^3}\phi_{u_\lambda}(x)u_\lambda v
dx\leqslant\int_{\mathbb{R}^{3}}u_\lambda^p v dx\ \text{ for all }
v\in H^1(\mathbb{R}^3), v\geqslant0,
\end{equation}
it follows that (\ref{eq:est1}) holds with $u=u_\lambda$ and
$\mu=1$. Hence, by Lemma \ref{Le:4.4.0}, we have
\begin{equation}\label{eq:th1.2.2}
\|u_\lambda\|_\infty\geq1 \text{ for all } \lambda>0.
\end{equation}
Meanwhile, for any nonnegative function $\varphi\in
C^\infty_0(\mathbb{R}^3)$ and any nonnegative piecewise smooth
function $h$ on $[0,+\infty)$, we see that $h(u_\lambda)\varphi\in
H^1(\mathbb{R}^3)$. Let $v=h(u_\lambda)\varphi$ in
(\ref{eq:th1.2.1}), it follows that (\ref{eq:4.3.1}) holds with
$u=u_\lambda$ and $N=3$. Hence, by Lemma \ref{Le:4.3}, we have
\begin{equation}\label{eq:th1.2.111}
\|u_\lambda\|_\infty\leqslant
C_1(1+\|u_\lambda\|_6^{C_2})\|u_\lambda\|_6.
\end{equation}
So, (\ref{eq:th1.2.2}) and (\ref{eq:th1.2.111}) imply that
$u_\lambda$ does not converge to 0 in $L^6(\mathbb{R}^3)$ as
$\lambda\rightarrow0$, then Theorems \ref{th:3.1} shows that there
exist nonnegative function $u\in E$ and $u\neq0$ such that,
\begin{equation} \nonumber
\int_{\mathbb{R}^{3}}\nabla u\nabla
\varphi+\frac{u\varphi}{|y|^{\alpha}}dx+\int_{\mathbb{R}^{3}}\phi_{u}(x)u\varphi
dx =\int_{\mathbb{R}^{3}}{u}^{p}\varphi dx,\ \text{ for all
}\varphi\in C_0^\infty(\mathbb{R}^3\setminus T).
 \end{equation}
Moreover, $w\in C^2(\mathbb{R}^3\setminus T)$. $\Box$

{\bf Acknowledgement:} This work was supported by
 NSFC(11071245,11171339 and 11126313). The second author thank also IMS-CUHK for its support during his visit in January 2012, where part of this work was carried out.


\begin{thebibliography}{10}

\bibitem{ASS}C.O.~Alves,~M.A.S.~Souto,~S.H.M.~Soares,~Schr\"{o}dinger-Poisson equations without Ambrosetti-Rabinowitz
condition,~J. Math. Anal. Appl.~377(2011)584--592.

\bibitem{AAmbroDRuiz-CCM}
A.~Ambrosetti,~D.~Ruiz,~Multiple bound states for the
{S}chr\"odinger-{P}oisson problem,~Commun. Contemp.
Math.~10(2008)391--404.

\bibitem{AAP}
A.~Azzollini,~P.~d'Avenia,~A.~Pomponio,~On the Schr\"{o}dinger
Maxwell equations under the effect of a general nonlinear term,~
Ann. Inst. H. Poincar\'{e} Anal. Non lin\'{e}aire~27(2010)779--791.

\bibitem{AAzzoAPom-Jmaa}
A.~Azzollini,~A.~Pomponio,~Ground state solutions for the nonlinear
{S}chr\"odinger-{M}axwell equations,~J. Math. Anal. Appl.~
345(2008)90--108.

\bibitem{BadMBenciVRolS-JEMS}
M.~Badiale,~V.~Benci,~S.~Rolando,~A nonlinear elliptic equation with
singular potential and applications to nonlinear field equations,~J.
Eur. Math. Soc.~9(2007)355--381.

\bibitem{MBadVBenSRol-reprent}
M.~Badiale,~V.~Benci,~S.~Rolando,~Three dimensional vortices in the
nonlinear wave equation,~Boll. Unione Mat. Ital.(9)~2(2009)105--134.

\bibitem{BadMGuiMRolS-AdvDE}
M.~Guida,~M.~Badiale,~S.~Rolando,~Elliptic equations with decaying
cylindrical potentials and power-type nonlinearities,~Adv. Differ.
Equ.~12(2007)1321--1362.


\bibitem{VBFD}V.~Benci,~D.~Fortunato,~Three dimensional vortices in
abelian gauge theories,~Nonlinear Anal.TMA~70(2009)4402-4421.

\bibitem{BFD1}
V.~Benci,~D.~Fortunato,~A minimization method and applications to
the study of solitons,~2011arXiv1111.1888.

\bibitem{BV}
V.~Benci,~N.~Visciglia,~Solitary waves with non-vanishing angular
momentum,~Adv. Nonlinear Stud.~3(2003)151-160.

\bibitem{GMColi-CAA}
G.M.~Coclite,~A multiplicity result for the nonlinear
{S}chr\"odinger-{M}axwell equations,~Commun. Appl. Anal.~
7(2003)417--423.


\bibitem{DT}
T.~D'Aprile,~On a class of solutions with non-vanishing angular
momentum for nonlinear Schr\"odinger equations, Differential
Integral Equations~16(2003)349-384.
\bibitem{TDApriDMugn-AdvNS}
T.~D'Aprile,~D.~Mugnai,~Non-existence results for the coupled
{K}lein-{G}ordon-{M}axwell equations,~Adv. Nonlinear
Stud.~4(2004)307--322.

\bibitem{TDApriDMugn-ProcRSES}
T.~D'Aprile,~D.~Mugnai,~Solitary waves for nonlinear
{K}lein-{G}ordon-{M}axwell and {S}chr\"odinger-{M}axwell
equations,~Proc. Roy. Soc. Edinburgh Sect. A~134(2004)893--906.

\bibitem{PDAven-AdvNS}
P.~D'Avenia,~Non-radially symmetric solutions of nonlinear
schr\"{o}dinger equation coupled with maxwell equations,~Adv.
Nonlinear Stud.~2(2002)177--192.

\bibitem{DavidGilbarg-Neil.S.Trudinger}
D.~Gilbarg,~N.S.~Trudinger,~Elliptic partial differential equations
of second order,~Classics in Mathematics. Springer-Verlag, Berlin,
2001, Reprint of the 1998 edition.

\bibitem{JZ-ActMS}
Y.S.~Jiang,~H.S.~Zhou,~Bound states for a stationary nonlinnear
Schr\"{o}dinger-Poisson systems with sign-changing potential in
$\mathbb{R}^3$,~Acta. Math. Sci.~{29B}(2009){1095--1104}.

\bibitem{JZ-reprint}
Y.S.~Jiang,~H.S.~Zhou,~Schr\"{o}dinger-Poisson system with steep
potential well,~J. Differential Equations~251(2011)582--608.

\bibitem{EHLiebMLoss-analysis}
E.H.~Lieb,~M.~Loss,~Analysis, volume~14 of Graduate Studies in
Mathematics, American Mathematical Society, Providence, RI, 1997.

\bibitem{Mug}
D.~Mugnai,~The Schr\"odinger-Poisson system with positive
potential,~Comm. Partial Differential Differential
Equations~36(2011)1099--1117.

\bibitem{DRuiz-JFA}
D.~Ruiz,~The {S}chr\"odinger-{P}oisson equation under the effect of
a nonlinear local term,~J. Funct. Anal.~237(2006)655--674.

\bibitem{SSol-AnnIHPA}
S.~Solimini,~A note on compactness-type properties with respect to
Lorentz norms of bounded subsets of a Sobolev space,~Ann. Inst. H.
Poincar\'{e} Anal. Non lin\'{e}aire~12(1995)319--337.

\bibitem{ZPWHSZhou-DisCDS}
Z.P. Wang,~H.S. Zhou,~Positive solution for a nonlinear stationary
{S}chr\"odinger-{P}oisson system in $\mathbb{R}^3$,~Discrete Contin.
Dyn. Syst.~18(2007)809--816.

\bibitem{MWil}
M.~Willem,~Minimax theorems,~Progress in Nonlinear Differential
Equations and their Applications, 24. Birkh\"auser Boston Inc.,
Boston, MA, 1996.

\bibitem{YZD}
M.B.~Yang,~Z.F.~Shen,~Y.H.~Ding,~Multiple semiclassical solutions
for the nonlinear Maxwell-Schr\"{o}dinger system,~Nonlinear Anal.
TMA~71(2009)730--739.

\bibitem{LgZhaoFkZhao-Jmaa}
L.G.~Zhao,~F.K.~Zhao,~On the existence of solutions for the
{S}chr\"odinger-{P}oisson equations,~J. Math. Anal. Appl.~
346(2008)155--169.
\end{thebibliography}

\end{document}